\RequirePackage{amsthm}
\documentclass[sn-mathphys-num]{sn-jnl}% Math and Physical Sciences Numbered Reference Style 
\usepackage{graphicx}%
\usepackage{multirow}%
\usepackage{amsmath,amssymb,amsfonts}%
\usepackage{lmodern}
\usepackage{mathrsfs}%
\usepackage[title]{appendix}%
\usepackage{xcolor}%
\usepackage{textcomp}%
\usepackage{manyfoot}%
\usepackage{booktabs}%
\usepackage{algorithm}%
\usepackage{algorithmicx}%
\usepackage{algpseudocode}%
\usepackage{listings}%
\usepackage{enumitem}
\usepackage{graphicx}
\usepackage{subfig}
\usepackage[misc]{ifsym}

    \makeatletter
\def\@fnsymbol#1{\ensuremath{\ifcase#1\or
   \mathsection\or \mathparagraph\or \|\or **\or \dagger\dagger
   \or \ddagger\ddagger \else\@ctrerr\fi}}
    \makeatother

\newcommand{\ra}[1]{\renewcommand{\arraystretch}{#1}}

\newcommand{\uv}[0]{\mathbf{u}}

\newcommand{\xv}[0]{\mathbf{x}}
\newcommand{\yv}[0]{\mathbf{y}}

\newcommand{\iprod}[2]{\left\langle {#1}, {#2} \right\rangle}
\newcommand{\norm}[1]{\left\lVert#1\right\rVert}

\newcommand{\dv}[0]{\mathbf{d}}
\newcommand{\ev}[0]{\mathbf{e}}
\newcommand{\av}[0]{\mathbf{a}}
\newcommand{\cv}[0]{\mathbf{c}}
\newcommand{\bv}[0]{\mathbf{b}}
\newcommand{\gv}[0]{\mathbf{g}}
\newcommand{\sv}[0]{\mathbf{s}}
\newcommand{\wv}[0]{\mathbf{w}}

\usepackage[normalem]{ulem}
\usepackage{cancel}
\newcommand{\stkout}[1]{\ifmmode\text{\sout{\ensuremath{#1}}}\else\sout{#1}\fi}

\newcommand{\coploweig}[0]{k}

\newcommand{\R}{\mathbb{R}}

\renewcommand{\SS}{\mathbb{S}}

\newcommand{\KK}{{K\!K}}
\newcommand{\CP}{\mathcal{CP}}
\newcommand{\COP}{\mathcal{COP}}
\renewcommand{\P}{\mathbb{P}}

\newcommand{\I}{\mathbb{I}}

\newcommand*\xbar[1]{%
  \hbox{%
    \vbox{%
      \hrule height 0.5pt % The actual bar
      \kern0.5ex%         % Distance between bar and symbol
      \hbox{%
        \kern-0.1em%      % Shortening on the left side
        \ensuremath{#1}%
        \kern-0.1em%      % Shortening on the right side
      }%
    }%
  }%
}

\theoremstyle{thmstyleone}%
\newtheorem{theorem}{Theorem}%  meant for continuous numbers
\newtheorem{proposition}[theorem]{Proposition}% 

\theoremstyle{thmstyletwo}%
\newtheorem{remark}{Remark}%

\theoremstyle{thmstylethree}%
\newtheorem{definition}{Definition}%
\newtheorem{lemma}{Lemma}

\begin{document}

\title[Article Title]{Sensitivity analysis for mixed binary quadratic programming\footnote{Parts of the paper have been published in proceedings of the 25th International Conference on Integer Programming and Combinatorial Optimization, IPCO 2024.}}

% \thanks{This is a titlenote}


\author[1]{\fnm{Diego} \sur{Cifuentes}}\email{dfc3@gatech.edu}
\equalcont{These authors contributed equally to this work.}

\author[1]{\fnm{Santanu} \sur{S. Dey}}\email{sdey30@gatech.edu}
\equalcont{These authors contributed equally to this work.}

\author*[1]{\fnm{Jingye} \sur{Xu}}\email{jxu673@gatech.edu}
\equalcont{These authors contributed equally to this work.}

\affil[1]{\orgdiv{H. Milton Stewart School of Industrial and Systems Engineering}, \orgname{Georgia Institute of Technology}, \orgaddress{\street{ 755 Ferst Dr NW}, \city{Atlanta}, \postcode{30332}, \state{GA}, \country{USA}}}


\abstract{We consider sensitivity analysis for Mixed Binary Quadratic Programs (MBQPs) with respect to changing right-hand-sides (rhs). We show that even if the optimal solution of a given MBQP is known, it is NP-hard to approximate the change in objective function value with respect to changes in rhs.
%it is NP-hard to achieve a two-sided $(a, b)$-approximation of the absolute value of the change in objective function value with respect to changes in rhs for any fixed value of $0 <a \leq 1\leq b$.
Next, we study algorithmic approaches to obtaining dual bounds for MBQP with changing rhs.
%Next, we examine the algorithmic challenges of producing good dual bounds for MBQPs with respect to changing rhs.
We leverage Burer's completely-positive (CPP) reformulation of MBQPs.
Its dual is an instance of co-positive programming (COP) and can be used to obtain sensitivity bounds.
%with respect to changing rhs.
We prove that strong duality between the CPP and COP problems holds if the feasible region is bounded or if the objective function is convex,
while the duality gap can be strictly positive if neither condition is met.
%When the feasible region is unbounded and the objective function is a non-convex quadratic, we show examples where there is a non-zero duality gap.
%We next show that the dimension of the set of optimal solutions of the COP dual is at least~$m$, where the original MBQP has $m$ constraints.
We also show that the COP dual has multiple optimal solutions,
and the choice of the dual solution affects the quality of the bounds with rhs changes.
We finally provide a method for finding good nearly optimal dual solutions,
and we present preliminary computational results on sensitivity analysis for MBQPs.}

\keywords{Sensitivity Analysis $\cdot$ Mixed Binary Quadratic Programming $\cdot$ Copositive programming $\cdot$ Duality Theory.}


\pacs[Mathematics Subject Classification]{90C11 $\cdot$ 90C31 $\cdot$ 90C46}

\maketitle

\section{Introduction.}
A mixed binary quadratic program (MBQP) has the form:
\begin{equation}\label{eq:MBQPeq}
\begin{aligned}
     z(\bv) := \min_{\xv \geq 0}\ & \xv^{\top} Q \xv + 2\cv^{\top} \xv \\
    \textup{s.t.}\ & \av_i^{\top} \xv = b_i,\ \forall i \in \{1, \dots, m\}\\
    &  x_j \in \{0,1\}, \forall j \in \mathcal{B},
    \end{aligned}
\end{equation}
where $Q$ is a symmetric matrix with rational entries of size $n \times n$, $\bv \in \mathbb{Q}^m$, $\cv\in \mathbb{Q}^n$, $\av_i \in \mathbb{Q}^n$ for all $i \in \{1, \dots, m\}$, and $\mathcal{B} \subseteq \{1, \dots, n\}$ is the set of variables restricted to be binary.
This is a very general optimization model that captures mixed binary linear programming~\cite{wolsey1999integer,conforti2014integer},
quadratic programming~\cite{bomze2002solving},
and several instances of mixed integer nonlinear programming models appearing in important application areas such as power systems~\cite{shahidehpour2003market}. 

Many practical optimization problems related to operational decision-making involve solving similar MBQP instances repeatedly. Moreover, they typically need to be solved within a short time window. This is because, unlike long-term planning problems, for such problems the exact problem data becomes available only a short time before a good solution is required to be implemented in practice. See, for examples, problems considered and discussed in~\cite{gu2022exploiting,johnson2020k,xavier2021learning}. In many of these applications, the constraint matrix remains the same, as these represent constraints related to some invariant physical resources, while the right-hand-side changes from instance to instance. 

The practical consideration discussed above motivates us to study sensitivity analysis of MBQPs with respect to changing right-hand-sides. {We are particularly interested in the change of the optimal objective value for (MBQP) given the fact that most modern solvers are able to find a good solution of (MBQP) in a reasonable time while 
finding it more challenging to provide matching dual bounds {(see preliminary experiments in appendix \ref{sec:dual_bounds_hard_to_find})}. 
Therefore, we aim to provide methods to infer high-quality dual bounds from sensitivity analysis of already-solved instances.} 
The pioneering results on sensitivity for integer programs (IPs) with changing right-hand-sides {were} obtained by Cook et al.~\cite{cook}. See~\cite{eisenbrand2019proximity,lee2020improving,del2022proximity,celaya2022improving,granot1990some}
for many advances in this line of research. However, these results yield trivial bounds in the case of binary variables since they rely on the infinity norm of integer constrained variables, which is a constant for all non-zero binary vectors.

We consider an alternative approach in this paper.
Specifically, we leverage Burer's completely-positive (CPP) reformulation of MBQP~\cite{burer}.
The advantage of the CPP reformulation is that, although still challenging and NP-hard in general to solve, it is a convex problem. Thus, one can examine its dual, which is an instance of copositive programming (COP). The optimal dual variables can provide bounds on $z(\bv)$, i.e., they allow to bound the optimal objective function of the MBQP as the right-hand-side changes.
Details of the CPP reformulation and the COP dual problem are presented in the next section.
This approach of using Burer's CPP reformulation of MBQPs~\cite{burer} to obtain shadow price information was first considered in~\cite{guocopositive} for the electricity market clearing problem.

%The rest of the organized as follows. 
{
The rest of the paper is organized as follows. section~\ref{sec:main} presents all our results.
Details of the proofs are presented in section~\ref{sec:hard} and section~\ref{sec:dual}. section~\ref{sec:conclude} provides conclusion and future avenues of research. 
}

An earlier version of this paper is published at the conference IPCO \cite{cifuentes2024sensitivity}. 
%The current paper is an improvement over \cite{cifuentes2024sensitivity}. 
%Although the main result (Theorem \ref{thm:allstrong}) is the same, we have reorganized the structure of section~\ref{sec:main} and have also added proofs to results in section~\ref{sec:main}, section~\ref{sec:hard} and section~\ref{sec:dual}. We have provided more detail on Vavasis' characterization of the optimal solution of quadratic programs \cite{VAVASIS199073} in section~\ref{sec:dual}. {We have also added more numerical experimental in section~\ref{sec:comp_body}}.

%\textcolor{red}{Add organization of the paper.}
%%%%%%%%%%%%%%%
%%%%%%%%%%%%%%%
%%%%%%%%%%%%%%%
%%%%%%%%%%%%%%%
%%%%%%%%%%%%%%%

\section{Main results.}\label{sec:main}
\paragraph{Notation.} Given a positive integer $n$, we let $[n]$ denote the set $\{1, \dots, n\}$.
For $u \in \mathbb{R}$, we denote its absolute value by $|u|$.
For a discrete set $\mathcal{B}$, we use $|\mathcal{B}|$ to denote its cardinality.
We let $\mathbb{S}^n$ to be the set of symmetric $n\times n$ matrices, and $\mathbb{S}^n_{+}$ to be the cone of $n\times n$ positive-semidefinite (PSD) matrices.
We denote a matrix $M$ being PSD by $M \succeq 0$ and denote $M$ not being a PSD matrix by $M \not \succeq 0$.
We let $\mathbb{S}^n_{p}$ be the set of symmetric $n\times n$ matrices with non-negative entries.
We let $\CP$ be the cone of completely positive matrices, i.e., $\CP = \{M \in \mathbb{S}^n\,|\, M = BB^{\top} \textup{ and }B \textup{ is a } m \times n \textup{ entry-wise nonnegative matrix for some integer }m\}$.
We let $\COP$ to be the cone of copositive matrices, i.e.,  $\COP = \{M \in \mathbb{S}^n\,|\,\xv^{\top} M \xv \geq 0, \forall \xv \geq 0\}$.
We use $\ev_i$ to denote the $i$-th standard basis vector.
{We say a matrix $M \in \COP$ is strictly copositive if $\xv^{\top} M \xv > 0, \forall \xv \geq 0$.} {For any symmetric matrix $M$, we use $\lambda_{\min}(M)$ to denote its smallest eigenvalue.}
{Throughout the paper, all vectors are assumed to be column vectors. Given two column vectors $\xv,\yv$, we use $[\xv;\yv]$ to denote the column vector by concatenating $\xv$ with $\yv$.}

\subsection{Complexity.}
We begin our study by formally establishing the difficulty of approximating $z(\bv + \Delta \bv)$ for varying $\Delta \bv$,
assuming that we know the exact value of $z(\bv)$. 
\begin{definition}
An algorithm is called ${(\Gamma_1,\Gamma_2)}$-approximation { of (\ref{eq:MBQPeq})} for some ${\Gamma_2} \geq 1 \geq {\Gamma_1} > 0$ if it takes $(A,\bv, \cv, Q, \mathcal{B}, z(\bv), \Delta \bv)$ as input,
where $A$, $\bv$, $\cv$, $Q$, $\mathcal{B}$ represents an instance of (\ref{eq:MBQPeq}), $z(\bv)$ is its optimal objective function value,
$\Delta \bv $ is the change in right-hand-side,
and it outputs a scalar $p$ satisfying
$${\Gamma_1} |\Delta z| \leq p \leq {\Gamma_2} |\Delta z|,$$ 
where $\Delta z = z(\bv) - z(\bv + \Delta \bv)$.
\end{definition}
We note that unlike the traditional definition of approximation for optimization, the two-sided bound is necessary. {Otherwise, if the lower bound {$\Gamma_1$} is not specified, an algorithm can ``cheat" by returning $p = 0$ to achieve $ p \leq  {\Gamma_2} |\Delta z|$. Similarly, if the upper bound {$\Gamma_2$} is not specified, an algorithm can ``cheat" by returning $p$ to be a very large number to achieve $ p \geq  {\Gamma_1} |\Delta z|$.}
% by returning either $p = 0$ or $p = \infty$ depending on whether $\alpha$ or $\beta$ is not specified. For example, to achieve $ p \leq  \beta |\Delta z|$, the algorithm can always return $p=0$. 

Our main result of this section is the following.
\begin{theorem}\label{thm:approx}
It is NP-hard to achieve $({\Gamma_1},{\Gamma_2})$-approximation { of (\ref{eq:MBQPeq})} {for general MBQPs for any ${\Gamma_2} \geq 1\geq {\Gamma_1} > 0$.}
\end{theorem}

Our proof of Theorem~\ref{thm:approx}, {presented in section~\ref{sec:hard},} is based on a reduction from the edge chromatic number problem,
using the fact that deciding whether the edge chromatic equals the max degree of a graph or one more than the max degree of a graph is NP-complete.
% Our proof is presented in Appendix~\ref{sec:hard}.

\subsection{Strong duality.}
We first present  results from~\cite{burer},
which are the starting point for our analysis. 
%Burer's reformulation holds for MBQPs when the constraints are in a specific form. Therefore, we begin by re-writing (\ref{eq:MBQP}) as follows
%\begin{eqnarray}\label{eq:MBQPeq}
%\begin{array}{rl}
%     z(\bv) & :=  \min  \xv^{\top} Q \xv + 2\cv^{\top} \xv \\
%    &  \ \ \ \ \ \av_i^{\top} \xv  = b_i,\ \forall i \in [m]\\
%    & \ \ \ \ \  x_j \in \{0,1\}, \forall j \in \mathcal{B}\\
%    & \ \ \ \ \  \xv\in \mathbb{R}^n_{+},
%    \end{array}
%\end{eqnarray}
%where we 
{The CPP reformulation of MBQP proposed in \cite{burer}}
% Burer's reformulation
makes the following assumption:
%assumes that the constraints $\xv \geq 0, \av_i^{\top} \xv = b_i,\ \forall i \in[m]$ implies that $0 \leq x_j \leq 1$ for all $j \in \mathcal{B}$.
\begin{equation}\tag{A}
\label{eq:assumption}
\xv \geq 0,\  \av_i^{\top} \xv = b_i,\ \forall i \in[m]
\quad\implies\quad
0 \leq x_j \leq 1 \text{ for all } j \in \mathcal{B}.
\end{equation}
As mentioned in~\cite{burer}, if $0 \leq x_j \leq 1$ for some $j \in \mathcal{B}$ is not implied, then we can explicitly add a constraint of the form $x_j + w_j = 1$ where $w_j \geq 0$ is a slack variable. Thus, this assumption is without any loss of generality.

Consider the following CPP problem:
\begin{equation}
\label{CP}
\begin{aligned}
z_{\CP}(\bv) :=  \min\ &\langle{C},{Y}\rangle \\
\text{s.t. } &\langle{T},{Y}\rangle = 1, \\
&\langle{A^{}_i},{Y}\rangle = 2 b_i,\forall i \in [m] \\
&\langle{AA^{}_i},{Y}\rangle = b_i^2 ,\forall i \in [m] \\
&\langle{N_j},{Y} \rangle= 0 , \forall j \in \mathcal{B} \\
&Y \in \CP,
\end{aligned}
\end{equation}
where $A_i \!=\! \begin{bmatrix} 0 & {\av}_i^{\!\top} \\ {\av}_i & 0  \end{bmatrix}$,
$AA_i \!=\! \begin{bmatrix} 0 & 0 \\ 0 & {\av}_i {\av}_i^{\!\top} \  \end{bmatrix}$,
$T \!=\! \begin{bmatrix}
    1 & 0 \\
    0 & 0
\end{bmatrix}$, 
$N_j \!=\! \begin{bmatrix}
    0 & -\ev_{j}^{\!\top} \\
    -\ev_{j} & 2 \ev_{j} \ev_{j}^{\!\top}
\end{bmatrix}$,
$C \!=\! \begin{bmatrix}
    0 & \cv^{\!\top}  \\
    \cv & Q  \\
\end{bmatrix}$.

Burer~\cite{burer} proves the following result.
\begin{theorem}[Burer's reformulation~\cite{burer}]\label{thm:Burer}
Given a feasible MBQP in the form (\ref{eq:MBQPeq}) satisfying assumption \eqref{eq:assumption},
%that $\xv \geq 0,\av_i^{\top} \xv = b_i,\ \forall i \in[m]$ implies that $0 \leq x_j \leq 1$ for all $j \in \mathcal{B}$,
then we have that $ z(\bv) = z_{\CP}(\bv)$.
\end{theorem}

Let us now consider the dual program to (\ref{CP})\footnote{For convenience, we have written the dual variables with 'negative sign'.}:
\begin{equation}
  \label{COP}
  \begin{aligned}
  z_{\COP}(\bv) := \textup{sup}\ &- \Bigl(\sum\limits_{i=1}^m 2b_i \alpha_i + b_i^2 \beta_i\Bigr)   - \theta \\
   \text{s.t. }\ &C + \sum\limits_{i=1}^m \Bigl(\alpha_i A_i + \beta_i AA_i\Bigr) + \Bigl(\sum\limits_{j \in \mathcal{B}} \gamma_j N_j \Bigr) + \theta T = M \\
   & M \in \COP. 
  \end{aligned}
\end{equation}

%Our goal, as outlined in the introduction, is the following:
Given an optimal solution to the dual (\ref{COP}), say $({\alpha^*}, {\beta^*}, {\gamma^*}, {\theta^*}, {M^*})$, and a perturbation to the right-hand-side of (\ref{eq:MBQPeq}) by $\Delta b\in \mathbb{R}^m$, we can obtain a lower bound to $z(\bv + \Delta b)$ as: 
\begin{eqnarray} \label{eq:weakdual}
z(\bv + \Delta \bv) \geq - \Bigl(\sum\limits_{i=1}^m 2(b_i + \Delta b_i) {\alpha^*_i} + (b_i + \Delta b_i)^2 {\beta^*_i}\Bigr)   - {\theta^*},
\end{eqnarray}
since this follows from weak duality. 

%In order for the the bound (\ref{eq:weakdual}) to be tighter, a nice property to have is that there is no duality gap between  (\ref{CP}) and  (\ref{COP}).
%This is the main topic of study in this section. 
If there is a positive duality gap between  (\ref{CP}) and  (\ref{COP}), then we do not expect the bound (\ref{eq:weakdual}) to be strong. Understanding when strong duality holds is the topic of this {{section}}. Our results of this {section} are aggregated in the next theorem:
\begin{theorem}[Strong duality]\label{thm:allstrong}
Consider a MBQP in the form (\ref{eq:MBQPeq}) satisfying the assumption \eqref{eq:assumption}.
%that  $\xv \geq 0, \av_i^{\top} \xv = b_i,\ \forall i \in [m]$ implies that $0 \leq x_j \leq 1$ for all $j \in \mathcal{B}$.
Let $\P = \{\xv :  \av_i^{\top} \xv  = b_i,\forall i \in [m], \xv \geq 0  \} \neq \emptyset$ denote the feasible region of (\ref{eq:MBQPeq}) that is assumed to be non-empty. Suppose $l$ is a finite lower bound on $z(\bv)$. Then:
\begin{enumerate}[label=(\alph*)]
    \item\label{p1} When $\P$ is bounded, there is a strictly copositive feasible solution of (\ref{COP}) which implies
    strong duality holds between (\ref{CP}) and (\ref{COP}) by the Slater condition.
    \item\label{p2} When $Q$ is PSD or $\P$ is bounded, for $\epsilon$ any arbitrarily small positive number, there is a closed-form formula to construct a feasible solution to (\ref{COP}) whose objective function value is $l-\epsilon$.
    In particular, when $l$ is the optimal value of (MBQP), this is a constructive proof that strong duality holds between (\ref{CP}) and (\ref{COP}).
    \item\label{p3} There exists examples where $Q$ is not PSD and $\P$ is unbounded, such that there is a positive duality gap between (\ref{CP}) and (\ref{COP}).
    \item\label{p4} The optimal solution of (\ref{COP}) is not attainable in general even if $\P$ is bounded and there is no duality gap.
\end{enumerate}
\end{theorem}
Note that part~\ref{p1} of Theorem~\ref{thm:allstrong} was shown in~\cite{brown2022copositive}, where the authors prove that strong duality holds between (\ref{CP}) and (\ref{COP}) in a non-constructive way when $\P$ is bounded. Their argument utilizes a recent result from~\cite{kim2021strong}. Also~\cite{LinderothRaghunathan} explores similar questions in a recent presentation.
To the best of our understanding, parts~\ref{p2}, \ref{p3}, and \ref{p4} of Theorem~\ref{thm:allstrong} were not known before.

We note that in part~\ref{p2} of Theorem~\ref{thm:allstrong} the promised closed form solutions can achieve additive $\epsilon$-optimal solutions for any $\epsilon >0$.
Is this an artifact of our proof technique? Clearly when $Q = 0$ and $\mathcal{B} = \emptyset$, the optimal solution of (\ref{COP}) can be achieved, as the dual optimal solution can be achieved by linear programming duality {(set $\gamma = 0, \beta = 0,\theta = 0$ and $\alpha \text{ equals linear program's dual variables}$)}.
Part~\ref{p4} shows that even if strong duality holds, the optimal solution is not attainable in (\ref{COP}) in general.
Therefore, part~\ref{p2} of Theorem~\ref{thm:allstrong} is the best we can hope for, as an optimal solution is not always achievable. Moreover, Part~\ref{p4} indicates that there is no Slater point in $(\ref{CP})$ in general when $\P$ is bounded. One can further show via simple examples that it is possible there is no Slater point in $(\ref{COP})$ when $\P$ is unbounded and $Q$ is PSD. 
% Combining those two examples, we can construct an example where strong duality holds between $(\ref{CP})$ and $(\ref{COP})$ while neither $(\ref{CP})$ nor $(\ref{COP})$ has a Slater point.
%%%%%%%%%%%%%%%
%%%%%%%%%%%%%%%

\subsubsection{Proof sketch for part~\ref{p1} of Theorem~\ref{thm:allstrong}: \ } 
{Let \begin{eqnarray}\label{eq:H}
H = T + \sum_{i = 1}^m AA_i.
\end{eqnarray}
For the proof of this part, when $\P$ is bounded, for a sufficiently large positive quantity $\hat{\lambda}$, 
we show that $\hat{M}:= C + \hat{\lambda}\cdot H$ is a strictly copositive matrix and a feasible solution of (\ref{COP}) where ${\alpha_i} = 0 \ \forall i \in [m]$,  ${\beta_i} = \hat{\lambda} \ \forall i \in [m]$, ${\gamma_j} = 0 \ \forall j \in \mathcal{B}$, ${\theta} = \hat{\lambda}$. Thus, this feasible solution is a Slater point leading to the required result. Details of our proof are presented in section~\ref{sec:boundedslater}.
}

%\sout{For the proof of this part, we show that $\hat{M}:= C + \hat{\lambda}\cdot H$ is a feasible solution of (\ref{COP}) and a strictly copositive matrix when $\P$ is bounded, where $\hat{\lambda}$ is sufficiently large positive quantity and}
%\begin{align*}
%\cancel{H = T + \sum_{i = 1}^m AA_i,}
%\end{align*}
%\sout{that is, $\hat{\alpha_i} = 0 \ \forall i \in [m]$,  $\hat{\beta_i} = \lambda \ \forall i \in [m]$, $\hat{\gamma_j} = 0 \ \forall j \in \mathcal{B}$, $\hat{\theta} = 0$.
%Thus $\hat{M}$ is a Slater point leading to the required result.}
% The full poof is in section~\ref{sec:boundedslater}.

%%%%%%%%%%%%%%%
%%%%%%%%%%%%%%%

\subsubsection{Proof sketch for part~\ref{p2} of Theorem~\ref{thm:allstrong}: \ } 

{Our proof of
part~\ref{p2} of Theorem~\ref{thm:allstrong} is long and therefore  presented in section~\ref{sec:constructsol}. Below we present some high-level details of the proof needed for further discussion in the paper.} The closed form solution promised in part~\ref{p2} of Theorem~\ref{thm:allstrong} is built using \emph{specific building blocks} or combinations of values for the variables $\alpha, \beta, \gamma, \theta$.
In particular we consider the following two building blocks:
%\begin{itemize}
\paragraph{(i) Building block 1.} For all $i \in [m]$, consider the following combination: $(\hat{\alpha_i} = -b_i, \hat{\beta_i} = 1, \hat{\theta} = b^2_i)$ and all other variables are zero; let the resulting matrix be: 
\begin{eqnarray}\label{eq:ri}
\KK_i = \sum\limits_{i=1}^m \left(\hat{\alpha_i} A_i + \hat{\beta_i} AA_i\right)   + \left(\sum\limits_{j \in \mathcal{B}} \hat{\gamma_j} N_j \right) + \hat{\theta} T = -b_i A_i + AA_i + b_i^2 T.
\end{eqnarray}
{
Note that $\KK_i$ is a $(n+1) \times (n+1)$ symmetric matrix, and given any $\yv := [t ; \xv] \in \R^{n+1}$, it is straightforward to verify that
\begin{align*}
   \yv^{\top} KK_i \yv = (b_i t - \av_i^\top \xv)^2,
\end{align*}
Therefore, one can interpret that $\KK_i$ is associated with the quadratic form obtained by homogenizing $(b_i - \av_i^\top \xv)^2$.}
%\sout{
%Note that $\KK_i$ is the the matrix associated to the quadratic form obtained by homogenizing $(b_i - \av_i^\top \xv)^2$.
%}
\paragraph{(ii)  Building block 2.} For all $j \in \mathcal{B}$, consider the following combination: $\tilde{\alpha_i} = -fb_i \ \forall i \in [m], \tilde{\beta_i} = f \ \forall i \in [m], \tilde{\theta} = (f\sum_{i = 1}^mb_i^2) + r,$ and $\tilde{\gamma_j} = -g$; let the resulting matrix be 
\begin{align*}
G_j(f, g, r) & = \sum\limits_{i=1}^m \Bigl(\tilde{\alpha_i} A_i + \tilde{\beta_i} AA_i \Bigr)   + \Bigl(\sum\limits_{j \in \mathcal{B}} \tilde{\gamma_j} N_j \Bigr) + \tilde{\theta} T \\ & = f\Bigl(\sum\limits_{i = 1}^m \KK_i\Bigr) - gN_j + rT,
\end{align*}
where $f, g, r$ are parameters. \\
%\end{itemize}

The closed-form solution that we construct for (\ref{COP})  is of the form:
\begin{eqnarray}\label{eq:dualsol}
U(f_1, f_2, g, r, \tau) = C + f_1\Bigl(\sum_{i = 1}^m \KK_i \Bigr)+ \sum_{j\in \mathcal{B}} G_j (f_2, g, r) + \tau H - l T, 
\end{eqnarray}
where $H$ is defined in (\ref{eq:H}). We specify values for the parameters $f_1, f_2, g, r, \tau$ such that the above matrix is copositive and has objective value of $l - \epsilon$. Our proof depends on the following Theorem \ref{thm_final_robust} stated below
that may be of independent interest. 

Consider the following perturbation of the original MBQP:
\begin{align}\tag{MBQP($\epsilon$)}
    \begin{array}{rl}
\zeta(\bv, \epsilon) := & \min\limits_{\xv,{\varepsilon}}  \xv^{\top} Q \xv +  2\cv^{\top} \xv \\
    &  \ \ \ \ \ \av_i^{\top} \xv  = b_i + \varepsilon^{(1)}_i,\forall i \in [m] \\
    & \ \ \ \ \  x_j + \varepsilon^{(2)}_j \in \{0,1\}, \forall j \in \mathcal{B} \\
    & \ \ \ \ \ {\varepsilon = [\varepsilon^{(1)} ; \varepsilon^{(2)}] \in \R^{m} \times \R^{|\mathcal{B}|},}  \\
    & \ \ \ \ \  \norm{\varepsilon^{(r)}}_{\infty} \leq \epsilon, \forall r \in \{1,2\},\\
    & \ \ \ \ \  \xv \geq 0.
    \end{array}
\end{align}
We prove the following result:
\begin{theorem}[Local stability]
\label{thm_final_robust} Let $l$ be a lower bound on $z(\bv)$, i.e., a lower bound on $\zeta(\bv, 0)$. When $Q \text{ is PSD}$ or $\P$ is bounded, there exists $t_1 > 0,t_2 \geq 0$ that depends on $A,\bv,\cv,Q,\mathcal{B}$ such that if $0 \leq \epsilon < t_1$, then $\zeta(\bv, \epsilon) \geq l - \epsilon t_2$.
\end{theorem}
{Our proof of Theorem~\ref{thm_final_robust} is presented in section~\ref{sec:localstable}.} We note that if we were only considering the case where $Q$ is PSD, then the above result could possibly be obtained using disjunctive arguments. Since we also allow for non-PSD $Q$ matrices (when $\P$ is bounded), our proof of Theorem~\ref{thm_final_robust} requires the use of result from~\cite{VAVASIS199073} characterizing the optimal solution of quadratic programs.
% , and is presented in Appendix~\ref{sec:localstable}.

% The details of our proof are presented in Appendix~\ref{sec:constructsol}. 

Note that by the construction of (\ref{CP}), the inequality $ z(\bv) \geq z_{\CP}(\bv)$ trivially holds. Part~\ref{p2} of Theorem~\ref{thm:allstrong} shows that $z_{\COP}(\bv) \geq z(\bv) - \epsilon$ for any positive $\epsilon$ when $Q$ is PSD or $\mathbb{P}$ is bounded. Since $z_{\CP}(\bv) \geq z_{\COP}(\bv) $ {is} a consequence of weak duality, we arrive at the following observation.

\begin{remark}[Alternative proof of Burer's Theorem]
The proof of part~\ref{p2} of Theorem~\ref{thm:allstrong} provides an alternative proof of Theorem~\ref{thm:Burer} in the case  when $Q$ is PSD or $\mathbb{P}$ is bounded. 
\end{remark}

%{Recent works consider solving (\ref{COP}) using cutting-plane techniques~\cite{ANSTREICHER2021218,badenbroek2022analytic,guocopositive,LinderothRaghunathan} as a way to solve  (\ref{COP}). 
We hope that the structure information provided by part~\ref{p2} of Theorem~\ref{thm:allstrong} could perhaps offer opportunities to design more efficient cutting-plane techniques.
%%%%%%%%%%%%%%%
%%%%%%%%%%%%%%%

\subsubsection{Proof for part~\ref{p3} of Theorem~\ref{thm:allstrong}: \ }
We will provide an example where $Q \not \succeq 0$ and $\mathbb{P}$ is unbounded,
such that (\ref{CP}) is feasible and has finite value while (\ref{COP}) is infeasible.
Consider the following instance:
%\begin{align*}
%    & \min x_1^2 - x_2^2 \\
%    & \text{s.t. } x_1 - x_2 = 0 \\
%    & \ \ \ \ \ \ x_1 \geq 0, x_2 \geq 0.
%\end{align*}
\begin{align*}
\min \{x_1^2 - x_2^2 \,|\, x_1 - x_2 = 0, x_1 \geq 0, x_2 \geq 0\}. \end{align*}
This problem is feasible and its optimal value is zero.
Hence, (\ref{CP}) is also feasible and has value zero by Theorem~\ref{thm:Burer}.
The COP dual is:
\begin{align*}
    \max &-\theta \\
    \text{s.t } & \begin{bmatrix}
        0 & 0 & 0 \\
        0 & 1 & 0 \\
        0 & 0 & -1 \\
    \end{bmatrix} + \theta \begin{bmatrix}
        1 & 0 & 0 \\
        0 & 0 & 0 \\
        0 & 0 & 0 \\
    \end{bmatrix} + \alpha \begin{bmatrix}
        0 & 1 & -1 \\
        1 & 0 & 0 \\
        -1 & 0 & 0  \\
    \end{bmatrix} + \beta \begin{bmatrix}
        0 & 0 & 0 \\
        0 & 1 & -1 \\
        0 & -1 & 1  \\
    \end{bmatrix} =: M \in \COP 
\end{align*}
We claim that the dual is infeasible.
Let $\yv = \begin{bmatrix}
    0 \\ 1 \\ 1+\epsilon
\end{bmatrix}$, where $\epsilon >0$.
Then
\begin{align*}
    \yv^{\top} M \yv & = 1 - (1+\epsilon)^2 + \beta(1+(1+\epsilon)^2-2(1+\epsilon))
     = - 2\epsilon  + (\beta - 1)\epsilon^2 
\end{align*}
When $\epsilon$ is small enough, $\yv^{\top} M \yv < 0$. This completes the proof. 
\begin{remark}[Local stability not satisfied when $Q\not \succeq0$ and $\mathbb{P}$ is unbounded]
It is instructive to see that the above example does not satisfy the local stability property.
Indeed, for any positive value of $\epsilon$, {we verify next that} $\zeta(0, \epsilon) = - \infty$, even though $\zeta(0,0) = 0$. {Observe that \begin{align*}
        \zeta(0, \epsilon) & = \min \{x_1^2 - x_2^2 \,|\, -\epsilon \leq x_1 - x_2 \leq \epsilon, x_1 \geq 0, x_2 \geq 0\} \\ & \leq x_1^2 - (x_1 +\epsilon)^2,\forall x_1 \geq 0 \text{ (choose $x_2 = x_1 + \epsilon$)} \\
        & = -2\epsilon x_1 + \epsilon^2,\forall x_1 \geq 0.
    \end{align*}
    As $x_1 \to \infty$, it follows that $\zeta(0, \epsilon)  \to - \infty$.}
Hence, the sufficient conditions for local stability Theorem~\ref{thm_final_robust} cannot be further relaxed.
\end{remark} 

%%%%%%%%%%%%%%%
%%%%%%%%%%%%%%%
\subsubsection{Part~\ref{p4} of Theorem~\ref{thm:allstrong}: \ } We prove this result for the copositive dual corresponding to the standard integer programming formulation for finding the stable set number of a clique with $6$ nodes, i.e., we consider the following MBQP:
%consider the maximum stable set problem for a graph $G = (V,E)$.
%Clearly, the optimal solution is just a single vertex.
%This can be written as the following (MBQP):
\begin{align*}
    & \min -2\sum_{j \in [6]} x_i \\
    & \text{s.t } x_u + x_v + s_{uv} = 1,\forall u < v, \ u, v \in [6]  \\ 
    &  \ \ \ \ \ \xv \in \{0,1\}^{6}, \sv \geq 0
\end{align*}
{Details of our proof is presented in section~\ref{sec:notatt}.}% \subsubsection{Part~\ref{p4} of Theorem~\ref{thm:allstrong}: \ } The proof is in Appendix~\ref{sec:notatt}. %Our poof of part~\ref{p4} of Theorem~\ref{thm:allstrong}  indicates that there is no slater point in $(\ref{CP})$ in general when $\P$ is bounded. One can also show that it is possible there is no slater point in $(\ref{COP})$ when $\P$ is unbounded and $Q$ is PSD. 
% Combining those two examples, we can construct an example where strong duality hold between $(\ref{CP})$ and $(\ref{COP})$ while neither $(\ref{CP})$ nor $(\ref{COP})$ has slater point.
%The detail is present in section~\ref{sec:notatt}.

% \begin{remark}
% \label{rem:no_slater_point}
% \end{remark}

%Finally, we note that in part~\ref{p2} of Theorem~\ref{thm:allstrong} the promised closed form solutions can achieve additive $\epsilon$-optimal solutions for any $\epsilon >0$. Is this an artifact of our proof technique? Clearly when $Q = 0$ and $\mathcal{B} = \emptyset$, the optimal solution of (\ref{COP}) can be achieved, as the dual optimal solution can be achieved by linear programming duality (set $\beta = 0$). However, it turns out that even if strong duality holds, the optimal solution is not attainable in (\ref{COP}) in general. Therefore, Part~\ref{p2} of Theorem~\ref{thm:allstrong} is the best we can hope for, as an optimal solution is not always achievable.
 
%%%%%%%%%%%%%%%
%%%%%%%%%%%%%%%

\subsection{How good is the closed form solution of Theorem~\ref{thm:allstrong} for sensitivity analysis?}\label{sec:multsol}

% Previous works consider solving (\ref{COP}) using cutting-plane techniques~\cite{ANSTREICHER2021218,badenbroek2022analytic,LinderothRaghunathan} as a way to solve the original MBQP or {solving (\ref{COP}) to design a pricing scheme in discrete energy markets \cite{guocopositive}}.
% However, in this paper we take a different perspective. We believe that with the success of modern state-of-the-art integer programming solvers, the original MBQP may be (in most cases) best solved directly using an integer programming solver.
The key attraction of Theorem~\ref{thm:allstrong} is to be able to build a closed-form solution (\ref{eq:dualsol})  of the dual~(\ref{COP}) using the optimal solution (or best known lower bound) of MBQP. 
{Therefore, one can directly start conducting }sensitivity analysis after solving the original MBQP and building the closed-form dual solution. 
{If the $\epsilon-$optimal dual solution is unique, then our constructed solution would be the only choice (for a fixed value of $\epsilon$). However, we will show next that the set of $\epsilon-$optimal dual solutions is never a singleton. Moreover, we shall show that the dual bound (when $\Delta b\neq 0$) implied by these different $\epsilon-$optimal dual solutions are not equal. In particular, the closed-form solutions will typically not provide the best dual bounds.}
%{In general, we wish to obtain a dual solution to provide a dual bound that is as large as possible.}

%However, conducting sensitivity analysis using dual solutions is challenging due to the presence of multiple $\epsilon$-optimal solution. 
First note that, given an $\epsilon$-optimal dual solution $(\alpha^*, \beta^*, \gamma^*, \theta^*)$ guaranteed by strong duality verified in Theorem~\ref{thm:allstrong}, we have that $z(\bv) = - \Bigl(\sum\limits_{i=1}^m 2b_i {\alpha^*_i} + b^2_i {\beta^*_i}\Bigr)   - {\theta^*} + \epsilon.$
% \begin{eqnarray*}
% z(\bv) = - \Bigl(\sum\limits_{i=1}^m 2b_i {\alpha^*_i} + b^2_i {\beta^*_i}\Bigr)   - {\theta^*} + \epsilon.
% \end{eqnarray*}
Subtracting the right-hand-side of the above from the right hand-side of (\ref{eq:weakdual}), we obtain that the predicted change in the objective function value using the dual solution $(\alpha^*, \beta^*, \gamma^*, \theta^*)$ when the right-hand-side changes form $\bv$ to $\bv + \Delta \bv$ is:
\begin{eqnarray}\label{eq:change_of_value}
 \textup{Predict}(\alpha^*, \beta^*, \gamma^*, \theta^*):= -\sum_{i=1}^m 2 \Delta b_i \alpha_i^* - \sum_{i=1}^m ((\Delta b_i)^2 + 2 b_i \Delta b_i)\beta_i^* - \epsilon. 
\end{eqnarray}
{Note that, $\textup{Predict}(\alpha^*, \beta^*, \gamma^*, \theta^*)$ is a lower bound on $z(b+ \Delta b) - z(b).$} Next consider the building block $\KK_i$ in \eqref{eq:ri},
used to construct the {closed form} solution in part~\ref{p2} of Theorem~\ref{thm:allstrong},
corresponding to $(\hat{\alpha_i} = -b_i, \hat{\beta_i} = 1, \hat{\theta} = b^2_i)$. 
% The objective function value of this block is 
% $- (2b_i \hat{\alpha_i} + b_i^2 \hat{\beta_i})   - \hat{\theta} = 0.$ Moreover, $\KK_i \succeq 0$ { and this implies $\KK_i \in \COP$. In this case, let $(\alpha^*, \beta^*, \gamma^*, \theta^*, M^*)$ be a feasible solution in of (\ref{COP}), it then follows that $M^* + \KK_i \in \COP$ as $\COP$ is a convex cone.} Thus, 
We arrive at the following observation {when we vary the contribution of $\KK_i$ in a solution}:
\begin{proposition}
\label{prop:recession1}
Let $\bv^i \in \mathbb{R}^m$ be the vector with $i^{\textup{th}}$ component equal to $b_i$ and zeros everywhere else. If $(\alpha^*, \beta^*, \gamma^*, \theta^*, M^*)$ is an $\epsilon$-optimal solution of (\ref{COP}), then $(\alpha^* - \bv^i, \beta^* + \ev_i, \gamma^*, \theta^* + b_i^2, M^* + \KK_i)$ is also an $\epsilon$-optimal solution of (\ref{COP}).
\end{proposition}

\begin{proof}
{
It is straightforward to verify that the objective value {contributed by} $\KK_i$ is $- (2b_i \cdot (-b_i) + b_i^2 )   - b_i^2 = 0$. Since the objective value of $(\alpha^* - \bv^i, \beta^* + \ev_i, \gamma^*, \theta^* + b_i^2, M^* + \KK_i)$ is the objective value of $(\alpha^*, \beta^*, \gamma^*, \theta^*, M^*)$ plus the objective value of $\KK_i$, this implies that $(\alpha^* - \bv^i, \beta^* + \ev_i, \gamma^*, \theta^* + b_i^2, M^* + \KK_i)$ and $(\alpha^*, \beta^*, \gamma^*, \theta^*, M^*)$ have the same objective value. Moreover, since $\KK_i \succeq 0$, this implies that $\KK_i \in \COP$ and so is $M^* + \KK_i$. Thus we  conclude that if $(\alpha^*, \beta^*, \gamma^*, \theta^*, M^*)$ is an $\epsilon$-optimal solution of (\ref{COP}), then $(\alpha^* - \bv^i, \beta^* + \ev_i, \gamma^*, \theta^* + b_i^2, M^* + \KK_i)$ is also an $\epsilon$-optimal solution of (\ref{COP}). 
}
\end{proof}

{Proposition~\ref{prop:recession1} shows that the set of $\epsilon$-optimal dual solution is always infinite. Now,} substituting $(\alpha^* - \bv^i, \beta^* + \ev_i, \gamma^*, \theta^* + b_i^2, M^* + \KK_i)$ in place of $(\alpha^*, \beta^*, \gamma^*, \theta^*, M^*)$ in (\ref{eq:change_of_value}) we obtain:
$$ \textup{Predict}(\alpha^* -\bv^i, \beta^* + \ev_i, \gamma^*, \theta^* + b_i^2) = \textup{Predict}(\alpha^*, \beta^*, \gamma^*, \theta^* )  - (\Delta b_i)^2.$$ 
% or equivalently, 
% $$ |\textup{Predict}(\alpha^* - \bv^i, \beta^* + \ev_i, \gamma^*, \theta^* + b_i^2)| = |\textup{Predict}(\alpha^*, \beta^*, \gamma^*, \theta^* )|  + (\Delta b_i)^2.$$
Thus, we arrive at the following conclusion:
\begin{remark} \label{rem:recession}
{Since $\textup{Predict}(\cdot)$ is a valid lower bound on $z(\bv + \Delta \bv) - z(\bv)$, we wish {for} this lower bound {to be}  as high as possible}. Therefore,
if $(\alpha^*, \beta^*, \gamma^*, \theta^*, M^*)$ is an $\epsilon$-optimal solutions of (\ref{COP}), then the lower bound obtained using the dual optimal solution $(\alpha^* - \bv^i, \beta^* + \ev_i, \gamma^*, \theta^* + b_i^2, M^* + \KK_i)$ for the right-hand-side vector $\bv + \Delta \bv^i$ is worse than that obtained by $(\alpha^*, \beta^*, \gamma^*, \theta^*, M^*)$.
\end{remark}

Therefore, in order to obtain the best possible sensitivity results,
we would like the contribution of $\KK_i$'s in the dual optimal matrix to be as small as possible.
The main role of $\KK_i$'s is to ensure that the constructed solution is in~$\COP$. 
However, as an artifact of our proof of part~\ref{p2} of Theorem~\ref{thm:allstrong},
the contribution of the $\KK_i$'s in the closed-form solution is much higher than what is really needed to ensure copositivity.
%In the proof of part~\ref{p2} of Theorem~\ref{thm:allstrong}, note that both the building blocks involve $\KK_i$'s (see (\ref{eq:ri}) and (\ref{eq:gj})).
%In fact, as an artifact of our proof of part~\ref{p2} of Theorem~\ref{thm:allstrong},
%the contribution of $\KK_i$ in the closed-form solution is much higher than
%a ``minimal" dual optimal solution with respect to contribution of $\KK_i$'s.
This fact was empirically verified by preliminary computations. 

By examining the structure of optimal solution $U(f_1, f_2, g, r, \tau)$ in (\ref{eq:dualsol})  and noting that the second building block $G_j(f, g,r)$ is a linear combination of $\KK_i$'s, $N_j$'s and $T$, we may try to find good dual solutions,
with small contribution of $\KK_i$'s
and fixed values of $\tau$ and $r$,
as follows:
\begin{equation}\label{eq:actual0}
\begin{aligned}
\min_{p, \gamma}\ & \sum_{i = 1}^m w_i p_i \\
\textup{s.t.}\ & C+ \sum_{i = 1}^m p_i \KK_i + \sum_{j \in \mathcal{B}}\gamma_j N_j + \tau H  - (l + r)T \in \COP,
\end{aligned}
\end{equation}
where $w_i$'s are some non-negative weights. 
{
The purpose of (\ref{eq:actual0}) is that we wish to find a dual solution with a small value of $p_i$. 
%Based on Remark \ref{rem:recession}, for any suboptimal dual solution in (\ref{eq:actual0}), there exists an optimal dual solution in (\ref{eq:actual0}) that provides better sensitivity analysis.
}
Note that part~\ref{p2} of Theorem~\ref{thm:allstrong} guarantees that the above problem finds an $\epsilon$-optimal ({the value of $\epsilon$ }depends on $\tau$ and $r$) dual solution. In our computations, we solved a variant of the above optimization problem.

%%%%%%%%%%%%%%%
%%%%%%%%%%%%%%%

\subsection{Preliminary computations.}
\label{sec:comp_body}
{In this section, we present experiments regarding the
dual bound of $z(\bv + \Delta \bv)$ for different classes of instances and choices of $\Delta \bv$.}
\subsubsection{Modifications to (\ref{eq:actual0}). \ }\label{sec:modify}
The following changes are made to improve the quality of the bound and the computational cost.
%We apply a few enhancements to (\ref{eq:actual0}). 

\paragraph{Linear penalty.}
Consider a new building block corresponding to $\alpha_i = -1, \theta = 2b_i$ and all other variables zero.
%This block has an objective function value $0$.
This block is associated to the homogenization of the linear function $(b_i - \av_i^\top \xv)$,
and it does not contribute to the objective function, just like $\KK_i$.
%One interpretation of this building block is that it is a linear penalty rather than a quadratic penalty like building block 1 and 2.
Setting, $K_i = 2b_iT - A_i$, we solve the following problem: 
\begin{equation}\label{eq:actual1}
\begin{aligned}
\min_{p, \gamma, \delta} \ & \sum_{i = 1}^m w_i^{(1)} p_i + w_i^{(2)} \delta_i  \\
\textup{s.t.}\ & C + \sum_{i = 1}^m p_i \KK_i + \sum_{i = 1}^m \delta_i K_i + \sum_{j \in \mathcal{B}}\gamma_j N_j  + \tau H  - (l + r)T \in \COP,
\end{aligned}
\end{equation}
Although including $K_i$ into the problem is not necessary, we have empirically observed that it leads to tighter sensitivity bounds due to the added degree of freedom. 
% The heuristic choice of $\wv^{(1)}$ and $\wv^{(2)}$ are presented in Appendix~\ref{sec:comp}.

\paragraph{Solving a restriction of (\ref{eq:actual1}).}
Solving a copositive program is challenging.
Therefore, we replaced the restriction of being in the copositive cone in (\ref{eq:actual1}) with a restriction of being in the $\mathbb{S}_{+} + \mathbb{S}_{P}$.
This leads to a semidefinite program, which can be solved in polynomial time.
However, the resulting problem can become infeasible.
To mitigate this problem, we consider two more changes:
\begin{enumerate}
\item Allowing non-optimal dual solutions: Instead of fixing $l$, we let $l$ become a variable. We also penalize finding a poor quality dual solution by changing the objective of 
(\ref{eq:actual1}) to: $\min_{p, \gamma, \delta,l}\ -l + \sum_{i = 1}^m w_i^{(1)} p_i + w_i^{(2)} \delta_i.$
% $$\min_{p, \gamma, \delta,l}\ -l + \sum_{i = 1}^m w_i^{(1)} p_i + w_i^{(2)} \delta_i.$$
In this way, we may increase the chances of finding a feasible solution, However the dual solution we find may be have lesser objective than the known optimal value of original MBQP. 
\item McCormick inequalities:
The $Y$ variable in (\ref{CP}) satisfies the following well-known McCormick inequalities:
\begin{align}
\label{eq:mccormick}
Y_{ij} \leq Y_{1,i}, \quad
Y_{ij} \leq  Y_{1,j}, \quad
Y_{ij} \geq Y_{1,i} + Y_{1,j} -1, \quad
Y_{ij} \geq 0.
\end{align}
We add new columns to (\ref{eq:actual1}) corresponding to these inequalities.
\end{enumerate}
{
In this case, for some given nonnegative vectors $\wv^{(1)},\wv^{(2)}$, our lower bounds are produced by optimizing the following modification of (\ref{eq:actual0}):
\begin{equation}
\label{eq:objectivefun}
    \begin{aligned}
\min_{p, \gamma, \delta,l,\iota} & \ -l + \sum_{i = 1}^m (w_i^{(1)} p_i + w_i^{(2)} \delta_i). \\
\textup{s.t.}\ & C + \sum_{i = 1}^m p_i \KK_i + \sum_{i = 1}^m \delta_i K_i + \sum_{j \in \mathcal{B}}\gamma_j N_j  + \tau H  - (l + r)T + \sum_{s} \iota_s F_s  \in \SS_+ + \SS_P, \\
&  \iota \geq 0,
    \end{aligned}
\end{equation}
where $\iota$ are dual variables associated with McCormick inequalities (\ref{eq:mccormick}) and each $F_s \in \SS^{n+1}$ corresponds to a McCormick inequality on $Y$.
}
% \begin{equation}\label{eq:objectivefun}
% \min_{p, \gamma, \delta,l}\ -l + \sum_{i = 1}^m (w_i^{(1)} p_i + w_i^{(2)} \delta_i).
% \end{equation}

\paragraph{Choice of $\wv^{(1)},\wv^{(2)}$.}

Though any vectors $\wv^{(1)},\wv^{(2)}$ allow us to derive lower bounds,
we can obtain better practical results by selecting them carefully.
Assume that the target range of $\Delta \bv_i$ is $\{0,1,\dots,\text{rg}_i\}$.
We select the vectors $\wv^{(1)},\wv^{(2)}$ as follows:
$$w_i^{(1)} = \dfrac{\sum_{\rho=0}^{\text{rg}_i}  \rho^2 }{\text{rg}_i+1} {= \dfrac{(2\text{rg}_i+1)(\text{rg}_i)}{6}}
\quad\text{ and }\quad
w_i^{(2)} =  \dfrac{\sum_{\rho=0}^{\text{rg}_i}  -2 \rho }{\text{rg}_i+1} = {-\text{rg}_i}.$$ 
The motivation behind such a choice is {that the primal problem is a minimization problem and we would like}
to maximize the average predicted lower bound over all $\Delta \bv$ in the target range,
as explained next.

Given optimal $l^*,p^*,\delta^*$, the predicted lower bound of $z(\bv + \Delta \bv)$ is 
\begin{align*}
 l^* - \sum\limits_{i=1}^m [p_i^{*} (\Delta b_i)^2 - 2 \delta_i^{*}  \Delta b_i].
\end{align*}
{
Note that the the average of $(\Delta b_i)^2$ over the target range is $(\sum_{\rho=0}^{\text{rg}_i}  \rho^2 )/(\text{rg}_i+1)$,
and the average of $-2 \Delta b_i$ is $(\sum_{\rho=0}^{\text{rg}_i} - 2 \rho)/(\text{rg}_i+1)$.
Hence, our choice of $\wv^{(1)}, \wv^{(2)}$ means that the objective value in \eqref{eq:objectivefun} corresponds to maximizing the average predicted lower bound over all the $\Delta \bv$ in the target range.}

\subsubsection{Preliminary experimental results.}

\paragraph{Instances.} In our preliminary experiments, we generate {five} classes of instances, which we refer to as (COMB), ({KBQP}), (SSLP), (SSQP), and {(PACK)}.
%For each class, we consider three density levels that control the sparsity of the constraints.

The first class of instances is a weighted stable set problem with a cardinality constraint:
\begin{equation}
\tag{COMB}
\min\bigl\{-\cv^{\top} \xv \mid \sum_{j = 1}^n x_j \leq p,\ x_i + x_j \leq 1\ \forall (i,j) \in E\bigr\}
\end{equation}
We generate random instances in the following way.
The underlying graph is a randomly generated bipartite graph $(V_1 \cup V_2, E) $ with $|V_1| = |V_2| = 10$ and each edge $(i,j) \in E$ is present with probability $d$ where $d \in \{0.3,0.5,0.7\}$. Each entry of $c$ is uniformly sampled from $\{0,\dots,10\}$. The right-hand-side of the cardinality constraint is $p=3$. 
Twenty instances were generated for each choice of $d$. For this class of instances, we performed sensitivity analysis with respect to the right-hand-side of the cardinality constraint, where we increased the value of $p$ by $\Delta p \in \{1, \dots, 10\}$.
%We refer this instances as (COMB) with density $d$ where $d \in \{0.3,0.5,0.7\}$.

{
The second class of instances is binary quadratic minimization with a cardinality constraint:
\begin{equation}
\tag{KBQP}
\min\bigl\{ \xv^{\top} Q \xv \mid \sum_{j = 1}^n x_j \leq p,\xv \in \{0,1\}^{n}\}
\end{equation}
We choose $p = 10$ and $Q$ from the BQP library generated by \cite{billionnet2007using} where $n=100$ and the density of $Q$ is $0.1$.
Ten instances were generated. For this class of instances, we performed sensitivity analysis with respect to the right-hand-side of the cardinality constraint, where we increased the value of $p$ by $\Delta p \in \{1, \dots, 10\}$.
}

The next class of instances contains continuous variables that are ``turned on or off" using binary variables. The instances have the following form:
\begin{equation}
\tag{SSLP}
\min\ \bigl\{
-2\cv_{x}^{\top}\! \xv + 2\cv_{y}^{\top}\! \yv
\mid
\av_i^{\top}\! \xv \!\leq\! b_i\ i \!\in\! [m],\ 
x_i \!\leq\! y_i\ i \!\in\! [n],\ 
x\!\geq\!0,\ y\!\in\! \{0,1\}^n
\bigr\}
\end{equation}
%\begin{equation}
%\tag{LPIND}
%\begin{aligned}
%\min\ &-2\cv_{\xv}^{\top} \xv + 2\cv_{\yv}^{\top} \yv + \xv^{\top} Q \xv \\
%\textup{s.t.}\ & \av_i^{\top} \xv \leq b_i,\ i \in [m],\\
%&x_i \leq y_i,\ i \in [n], \\
%& x\geq0,\ y\in \{0,1\}^n.
%\end{aligned}
%\end{equation}
We generate instances in the following way.
We set $n=20$ and $m = 5$.
Each entry of $\cv_{x}$ is uniformly sampled from $\{0,\dots,10\}$ and $\cv_{y} = (3,\dots,3)$ is a constant vector.
Each entry of $\av_i$ is uniformly sampled from $\{0,\dots,10\}$ and then each entry of $\av_i$ is zeroed out with probability $d \in \{0.3, 0.5, 0.7\}$.
Finally, $b_i = \lfloor \frac{1}{2} \av_i^{\top} \ev \rfloor$ for all $i \in [m]$.
Twenty instances were generated for each probability. For this class of instances, we focus on sensitivity with-respect to right-hand-side of $\av_i^{\top}x \leq b_i,\forall i \in [m]$ when $\Delta \bv \in \{0,1,2,3\}^m$.
%We refer this instances as (LPIND) with density $d$ where $d \in \{0.3,0.5,0.7\}$.

The next class of instances is similar,
except that the objective is quadratic:
\begin{equation}
\tag{SSQP}
\min\ \bigl\{
-2\cv_{x}^{\top}\! \xv + 2\cv_{y}^{\top}\! \yv + \xv^{\top}\! Q \xv
\mid
\av_i^{\top}\! \xv \!\leq\! b_i\ i \!\in\! [m],\ 
x_i \!\leq\! y_i\ i \!\in\! [n],\ 
x\!\geq\!0,\ y\!\in\! \{0,1\}^n
\bigr\}
\end{equation}
We generate $\cv_{x}$, $\cv_{y}$, $\av_i$,$\bv$ as before.
The matrix $Q$ is randomly generated such that $Q = \sum\limits_{i \in \{1,2\}} \uv_i \uv_i^{\top}$ where each entry of $\uv_i$ is uniformly sampled from $\{-1,0,1\}$.

{
The last class of instance is a binary packing problem:
\begin{equation}
\label{eq:PACK}
\tag{PACK}
\min\ \bigl\{
-\cv^{\top}\! \xv
\mid
\av_i^{\top}\! \xv \!\leq\! b_i\ i \!\in\! [m],\ x\!\in\! \{0,1\}^n
\bigr\}
\end{equation}
%\begin{equation}
%\tag{LPIND}
%\begin{aligned}
%\min\ &-2\cv_{\xv}^{\top} \xv + 2\cv_{\yv}^{\top} \yv + \xv^{\top} Q \xv \\
%\textup{s.t.}\ & \av_i^{\top} \xv \leq b_i,\ i \in [m],\\
%&x_i \leq y_i,\ i \in [n], \\
%& x\geq0,\ y\in \{0,1\}^n.
%\end{aligned}
%\end{equation}
We generate instances in the following way.
We set $n=100$ and $m = 10$.
Each entry of $\cv$ is uniformly sampled from $\{0,\dots,3\}$ and
each entry of $\av_i$ is uniformly sampled from $\{0,\dots,3\}$ and then each entry of $\av_i$ is zeroed out with probability $d \in \{0.3, 0.5, 0.7\}$.
We choose $b_i = \lfloor \frac{1}{2} \av_i^{\top} \ev \rfloor$ for all $i \in [m]$. We solve (\ref{eq:PACK}) to get optimal solution $\xv^*$ and then update $\bv = A \xv^* - \ev$. We updated $b$, in this way because otherwise we empirically observed that (\ref{eq:PACK}) with randomly generated data is likely to be insensitive to small changes of the right-hand-side.
Ten instances were generated for each probability. For this class of instances, we focus on sensitivity with-respect to right-hand-side of $\av_i^{\top}x \leq b_i,\forall i \in [m]$ when $\Delta \bv \in \{0,1,2,3\}^m$. However, due to the curse of dimensionality, we cannot examine all instance with changed right-hand-side.  Therefore, we uniformly sample 2000 possible choice of $\Delta \bv \in \{0,1,2,3\}^m$ and compare sensitivity analysis with different methods.}

% We consider two experiments. For the first experiments, each $b_i$ is half of the sum of the left-hand-side entries, $p = 3$ and we focus on sensitivity of $\ev^{\top} \yv \leq p$ when $\Delta p \in \{1,\dots,10\}$. For the second experiment, each $b_i$ is a third of the sum of the left-hand-side entries, $p = 10$ and we focus on sensitivity with-respect to right-hand-side of $\av_i^{\top} \leq b_i,\forall i \in [m]$ when $\Delta \bv \in \{0,1,2,3\}^m$. \textcolor{red}[To discuss]
% \paragraph{Software and hardware:}

% Our experiments are implemented in Julia 1.9, relying
% on Gurobi version 9.0.2 and Mosek 10.1 as the solvers.
% We solve on a Windows PC with 12th Gen Intel(R) Core(TM) i7 processors and 16 RAM.
\paragraph{Experiments conducted.}
Our experiments are implemented in Julia 1.9, relying
on Gurobi version 9.0.2 and Mosek 10.1 as the solvers.
We solve on a Windows PC with 12th Gen Intel(R) Core(TM) i7 processors and 16 RAM.
We compare our method with other known methods.
Those methods consider certain convex relaxation of (MBQP) and obtain dual variables of constraints to conduct sensitivity analysis via weak duality.
In this case, we consider three convex relaxations,
which we call Shor1, Shor2, Cont. {Shor1 and Shor2 are standard Shor's semidefinite relaxations of MBQP \cite{shor1987quadratic} with additional non-negativity constraint and standard McCormick inequalities \cite{mccormick1976computability}. Cont is a naive convex relaxation by dropping integral constraints. We present these relaxations in detail below.}
%first order Shor relaxation, second order Shor relaxation (from dual side, it can be considered as a naive restriction of (\ref{COP})) and convex quadratic relaxation.

{The first relaxation of (MBQP), referred to as Shor1 is the following:}
%First, we consider the the Shor relaxation of (MBQP):
%     \[
%     \tag{Shor1}
% \min_{Y \in \SS_+ \cap \SS_P} \left\{ \iprod{C}{Y}  \;\middle\vert\;
%    \begin{array}{@{}l@{}}
% \iprod{T}{Y} = 1, \iprod{A_i}{Y} = 2 b_i, \forall i \in [m],\\
% \iprod{N_j}{Y} = 0,\forall j \in \mathcal{B}
%    \end{array}
% \right\}.
% \]
    \[
    \tag{Shor1}
\min_{Y \in \SS_+ \cap \SS_P} \left\{ \iprod{C}{Y}  \;\middle\vert\;
   \begin{array}{@{}l@{}}
\iprod{T}{Y} = 1, \iprod{N_j}{Y} = 0,\forall j \in \mathcal{B} \\
\iprod{A_i}{Y} = 2 b_i, \forall i \in [m] \\
\text{McCormick inequalities (\ref{eq:mccormick})} 
   \end{array}
\right\}.
\]
% \begin{equation}
% \tag{Shor1}
% \min\ \bigl\{ \iprod{C}{Y} | \iprod{T}{Y} = 1, \iprod{A_i}{Y} = 2 b_i, \forall i \in [m], \\ \iprod{N_j}{Y} = 0,\forall j \in \mathcal{B},Y \geq 0, Y \succeq 0. \bigr\} 
% \end{equation}
% \begin{equation}
% \tag{Shor1}
% \begin{aligned}
% \min & \iprod{C}{Y} \\
%  \textup{s.t. } &  \iprod{T}{Y} = 1, \\
% & \iprod{A_i}{Y} = 2 b_i, \forall i \in [m],  \\
% & \iprod{N_j}{Y} = 0,\forall j \in \mathcal{B}, \\
% & Y \geq 0, Y \succeq 0.
% \end{aligned}
% \end{equation}
Next, we consider the relaxation of the problem obtained by augmenting Shor1 with redundant quadratic constraints and McCormick inequalities:
    \[
    \tag{Shor2}
\min_{Y \in \SS_+ \cap \SS_P} \left\{ \iprod{C}{Y}  \;\middle\vert\;
   \begin{array}{@{}l@{}}
\iprod{T}{Y} = 1, \iprod{N_j}{Y} = 0,\forall j \in \mathcal{B} \\
\iprod{A_i}{Y} = 2 b_i, \iprod{AA_i}{Y} = b_i^2, \forall i \in [m] \\
\text{McCormick inequalities (\ref{eq:mccormick})} 
   \end{array}
\right\}.
\]
% \begin{equation}
% \tag{Shor2}
% \begin{aligned}
% \min\ & \iprod{C}{Y} \\
% \textup{s.t.}\ & \iprod{T}{Y} = 1, \\
% & \iprod{A_i}{Y} = 2 b_i, \forall i \in [m],   \\
% & \iprod{AA_i}{Y} = b_i^2, \forall i \in [m],   \\
% & \iprod{N_j}{Y} = 0,\forall j \in \mathcal{B}, \\
% & Y \geq 0, Y \succeq 0.
% \end{aligned}
% \end{equation}
% \textcolor{red}{McCormick?}
%We solve the above semidefinite program using

Finally, and assuming that $Q$ is PSD,
we obtain a convex relaxation simply by relaxing the binary variables to be continuous variables in~$[0,1]$:
    \[
    \tag{Cont}
\min_{\xv \geq 0} \left\{ \xv^{\top} Q \xv + 2\cv^{\top} \xv  \;\middle\vert\;
   \begin{array}{@{}l@{}}
\av_i^{\top} \xv = b_i,\ \forall i \in [m],x_j \in [0,1] \text{ and } \forall j \in \mathcal{B}
   \end{array}
\right\}.
\]
% \begin{equation}
% \tag{Cont}
% \begin{aligned}
% \min\ & \xv^{\top} Q \xv + 2\cv^{\top} \xv \\
% \textup{s.t.}\ & \av_i^{\top} \xv = b_i,\ \forall i \in [m],\\
% &  x_j \in [0,1], \forall j \in \mathcal{B},\\
% &  \xv\in \mathbb{R}^n_{+},
% \end{aligned}
% \end{equation}

{In section~\ref{sec:modify}, we discussed the choice of $w^{(1)}$ and $w^{(2)}$. In order to test the performance of this choice, we also tested other values of  $w^{(1)} $ and $w^{(2)} $ for (COMB) instances.} {In particular, we considered $w^{(1)} = (\sigma)^2$ 
 and $w^{(2)}  = -2(\sigma)$, where $\sigma \in \{2, 4, 6, 8, 10\}$. }

The relaxations Shor1 and Shor2 are solved using Mosek,
while the relaxation Cont is solved using Gurobi.
Notice that these problems may have multiple optimal solutions, so different solvers might lead to different solutions.

%Given $\Delta \bv$, our method along with those three methods will provide a prediction (lower bound) on $z(\bv + \Delta \bv)$.
We choose relaxation Shor1 as a baseline and measure the goodness of those predictions by \textit{relative gap}.
Given the rhs change $\Delta \bv$, the ground-truth $z(\bv + \Delta \bv)$, the prediction $p_1$ by Shor1 and the prediction $p_2$ by some method,
%the relative gap for this method at this particular rhs change $\Delta \bv$ is
then
$$\text{relative gap} = \frac{z(\bv + \Delta \bv) - p_2}{z(\bv + \Delta \bv) - p_1}.$$
This is always a non-negative number and a smaller relative gap indicates a better performance of the given method.

\begin{figure}[h!]%
    \centering
    \subfloat[\centering  a (COMB) instance with $d=0.3$]{{\includegraphics[width=5.5cm]{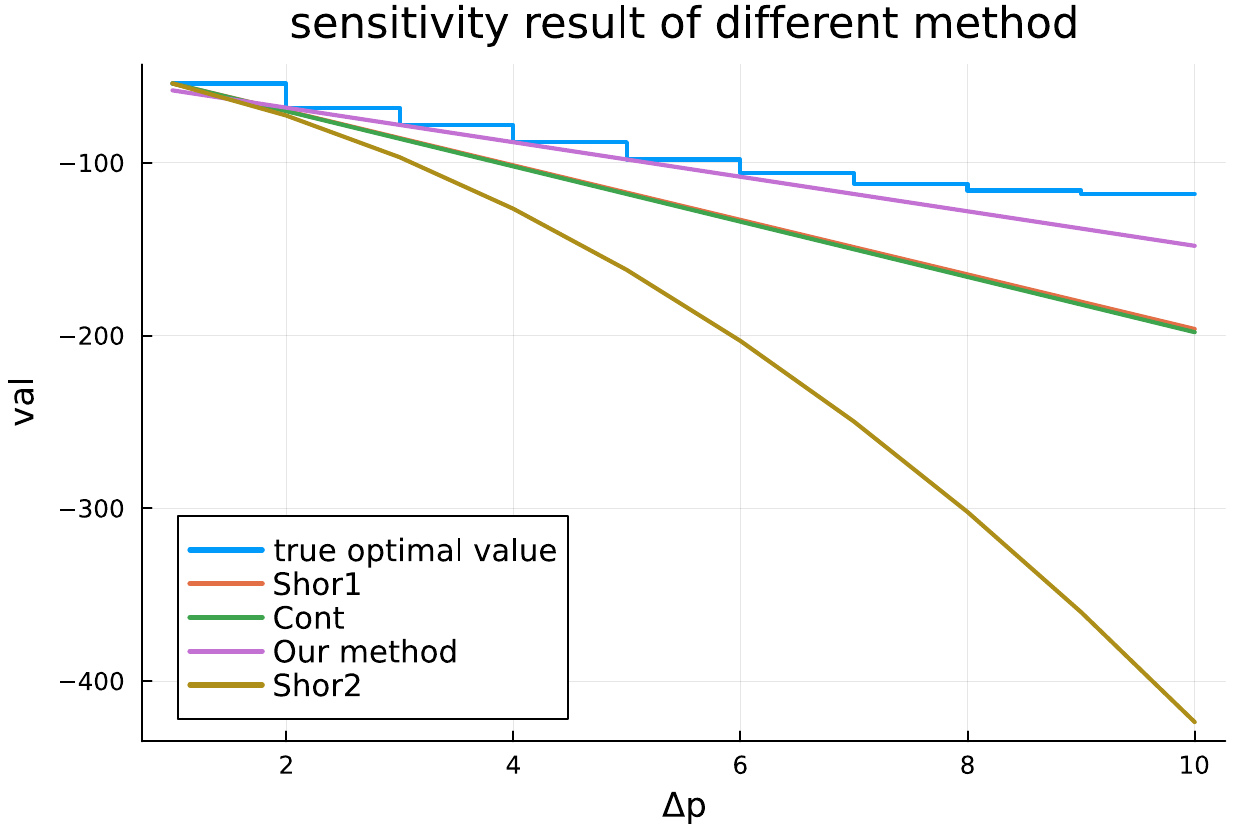} }}%
    \qquad
    \subfloat[\centering   a (COMB) instance with $d=0.3$]{{\includegraphics[width=5.5cm]{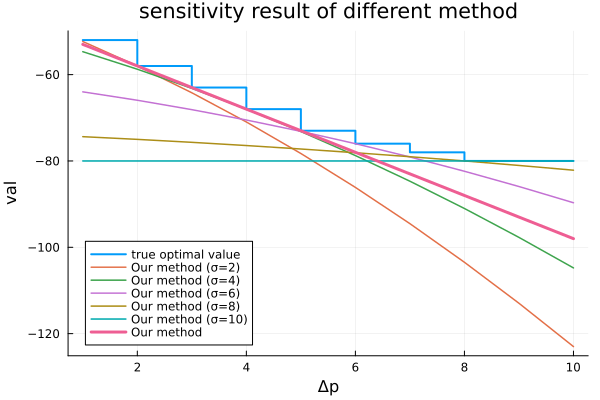} }}%
    \caption{The x-axis corresponds to $\Delta p$, and the y-axis corresponds to the optimal value of the new program or different predicted values of different methods. {Figure (a) plots Shor1, Shor2, Our method, Cont for a (COMB) instance with $d=0.3$.
    Figure (b) plots Our method, Our method($\sigma = 2$), Our method($\sigma = 4$), Our method($\sigma = 6$), Our method($\sigma = 8$) for a (COMB) instance with $d=0.3$.}
    }%
    \label{fig:parabla}
\end{figure}

\paragraph{Results and discussion.}
Figure~\ref{fig:parabla} shows an example of bounds obtained using different methods. 
{Tables~\ref{tab:comb}, \ref{tab:SQP}, \ref{tab:lpind},  \ref{tab:qpind}, and \ref{tab:PACK} }summarize the relative gaps obtained in the experiments mentioned above.
{See Appendix~\ref{Appdix:comp} for more detailed version of the results where the relative gaps are presented separately for each instance type and density of data.}

% Appendix~\ref{sec:comp} provides a more detailed information, including relative gaps per density.

{\footnotesize
\begin{table}[h!]\centering
`\caption{Average relative gap (COMB) -- all densities}
\label{tab:comb}
% \ra{0.5}
\setlength{\tabcolsep}{4pt}
\begin{tabular}{@{}cccccccccccccc@{}}\toprule 
$\Delta k$
& 1 & 2 & 3 & 4 & 5 & 6 & 7 & 8 & 9 & 10 & \phantom{a} & time(s)\\
\midrule
Shor1 & 1  & 1 & 1 & 1 & 1 & 1 & 1 & 1 & 1 & 1 & & 3.13 \\
Shor2   & 1.63&	3.71&	5.53	&6.90&	8.09&	6.75&	6.25&	6.08&	5.94	&5.84	&&4.15 \\
Our method &  1.28&	0.35&	0.06&	0.00	&0.07&	0.21	&0.29&	0.37	&0.43&	0.48&	&4.59\\
Our method ($\sigma =2$) &  0.70	&0.01	&0.32&	0.90	&1.54&	1.70	&1.86&	2.02	&2.14&	2.23&&4.68\\
 Our method ($\sigma =4$) &  4.76	&1.23&	0.21&	0.00	&0.1	&0.40	&0.58	&0.75&	0.89&	0.98&&	4.96\\
 Our method ($\sigma =6$) &  14.32&  5.27& 	2.11& 	0.66& 	0.10& 	0.00	& 0.047	& 0.13	& 0.23	& 0.33&& 	4.23
\\
 Our method ($\sigma =8$) &  23.93	&10.14	&4.68	&2.09&	0.84&	0.20&	0.03&	0.00&	0.02	&0.07&&	4.81\\
 Our method ($\sigma =10$) &  35.15	& 15.87&	8.40&	4.35	&2.24&	0.803&	0.30	&0.09&	0.02&	0.00&&4.40
\\
Cont   &  0.90	&1.01&	1.06	&1.03	&1.04	&1.04&	1.03	&1.04	&1.03	&1.02&&	0.00 \\
\bottomrule
\end{tabular}
\end{table}
}

{
{
{\footnotesize
\begin{table}[h!]\centering
\caption{Average relative gap (SBQP)}
\label{tab:SQP}
% \ra{0.5}
\setlength{\tabcolsep}{4pt}
\begin{tabular}{@{}cccccccccccccc@{}}\toprule 
$\Delta k$
& 1 & 2 & 3 & 4 & 5 & 6 & 7 & 8 & 9 & 10 & \phantom{a} & time(s)\\
\midrule
Shor1 & 1  & 1 & 1 & 1 & 1 & 1 & 1 & 1 & 1 & 1 & & 1313 \\
Shor2   & 4.80 & 2.52 & 2.17 & 2.13 & 2.06 & 2.06 & 2.12 & 2.15 & 2.15 & 2.13
&& 1154 \\
Our method &  0.13 & 0.02 & 0.20 & 0.36 & 0.48 & 0.55 & 0.60 & 0.63 & 0.67 & 0.70  && 1845 \\
\bottomrule
\end{tabular}
\end{table}
}
}
}

\begin{table}[h!]
    \parbox{.45\linewidth}{\footnotesize
\caption{Average relative gap for (SSLP) -- all densities}
\label{tab:lpind}
% \ra{0.5}
\setlength{\tabcolsep}{3pt}
\centering
 \begin{tabular}{ccccc} 
 \toprule
$\norm{\Delta \bv}_{\infty}$ &$\leq 1$ & $\leq 2$ & $\leq 3$ & time(s)  \\  
  \midrule
 Shor1 & 1 & 1 & 1 &3.63 \\
 Shor2 &  1.20  &  1.48 &  1.64 & 7.21 \\
 Our method  & 0.59  &  0.55 &  0.60 &	 5.82 \\
 Cont & 1.00 & 1.00 & 1.00 & 0.00 \\ [1ex] 
 \bottomrule
 \end{tabular}

}
\hfill
\parbox{.45\linewidth}{\footnotesize

\caption{Average relative gap for (SSQP) -- all densities}
\label{tab:qpind}
% \ra{0.5}
\setlength{\tabcolsep}{3pt}
\centering
 \begin{tabular}{ccccc} 
 \toprule
$\norm{\Delta \bv}_{\infty}$ &$\leq 1$ & $\leq 2$ & $\leq 3$ & time(s) \\  
  \midrule
 Shor1 & 1 & 1 & 1 & 3.58 \\
 Shor2 &  1.24  &  1.39 &  1.54	 & 7.08 \\
 Our method  & 0.52   & 0.48 &   0.53 & 5.90 \\
 Cont & 1.00 & 1.00 & 1.00 & 0.00 \\ [1ex] 
 \bottomrule
 \end{tabular}
}
\end{table}

\begin{table}[h!]
    {\footnotesize
\caption{Average relative gap of 2000 samples for (PACK) -- all densities}
\label{tab:PACK}
% \ra{0.5}
\setlength{\tabcolsep}{20pt}
\centering
 \begin{tabular}{ccc} 
 \toprule
$\norm{\Delta \bv}_{\infty}$ & $\leq 3$ & time(s)  \\  
  \midrule
 Shor1 & 1  &890 \\
 Shor2 &  3.53   &1139 \\
 Our method  & 0.60     &	 1001 \\
 Cont & 0.98 &  0.00 \\ [1ex] 
 \bottomrule
 \end{tabular}

}
\end{table}
% {\footnotesize
% \begin{table}[h!]
% \caption{Average relative gap for (SSLP) -- all densities}
% \label{tab:lpind}
% % \ra{0.5}
% \setlength{\tabcolsep}{6pt}
% \centering
%  \begin{tabular}{ccccc} 
%  \toprule
%  %\multicolumn{4}{c}{\text{Experimental summary: (LPIND) - average relative gap - all densities}} \\ 
%  %\midrule
% $\norm{\Delta \bv}_{\infty}$ &$\leq 1$ & $\leq 2$ & $\leq 3$ & avg time(s)  \\  
%   \midrule
%  Shor1 & 1 & 1 & 1 &3.63 \\
%  Shor2 &  1.20  &  1.48 &  1.64 & 7.21 \\
%  our method  & 0.59  &  0.55 &  0.60 &	 5.82 \\
%  Cont & 1.00 & 1.00 & 1.00 & - \\ [1ex] 
%  \bottomrule
%  \end{tabular}
% \end{table}
% }
% \vfill
% {\footnotesize
% \begin{table}[h!]
% \caption{Average relative gap for (SSQP) -- all densities}
% \label{tab:qpind}
% % \ra{0.5}
% \setlength{\tabcolsep}{8pt}
% \centering
%  \begin{tabular}{ccccc} 
%  \toprule
%  %\multicolumn{4}{c}{\text{Experimental eragsummary: (QPIND) - average relative gap - all densities}} \\ 
%  %\midrule
% $\norm{\Delta \bv}_{\infty}$ &$\leq 1$ & $\leq 2$ & $\leq 3$ & avg time(s) \\  
%   \midrule
%  Shor1 & 1 & 1 & 1 & 3.58 \\
%  Shor2 &  1.24  &  1.39 &  1.54	 & 7.08 \\
%  our method  & 0.52   & 0.48 &   0.53 & 5.90 \\
%  Cont & 1.00 & 1.00 & 1.00 & - \\ [1ex] 
%  \bottomrule
%  \end{tabular}
% \end{table}
% }
We observe that our method provides the tightest sensitivity bounds in all cases. {Also observe that our choice of $w^{(1)}$ and $w^{(2)}$ performs quite well across different values of $\Delta b$ as seen in Figure~\ref{fig:parabla}.}
Also note that method Shor2 provides the worst bounds. This is easier to observe in Figure \ref{fig:parabla}.
This is interesting, because the SDP from Shor2 is quite similar to Burer's formulation (with additional McCormick inequalities).
This discrepancy is most likely due to the fact that these problems have multiple optimal solutions - similar to the discussion in section~\ref{sec:multsol}.
The naive Shor2 approach finds an optimal dual which does not give good bounds after $\bv$ is perturbed.
On the other hand, our method attempts to find a good dual solution (with respect to producing good bounds for changing rhs) inside the $\epsilon$-optimal face of the dual.

{Although computing dual solution takes a considerable time, once such dual solution is available, one can compute the dual bound on $z(\bv + \Delta \bv)$ by simple arithmetic calculations, while computing $z(\bv + \Delta \bv)$ directly requires solving a MBQP. Therefore, the former approach is more scalable with the increase in the dimension of $\Delta \bv$.
For example, for the (\ref{eq:PACK}) instances above, computing the dual bound of $z(\bv + \Delta \bv)$ by dual solution for all $\Delta \bv \in \{0,1,2,3\}^{10}$ takes around 10 seconds while computing $z(\bv + \Delta \bv)$ by solvers is expected to take more than one day as the number of choice of $\Delta \bv$ is $4^{10}$.}


% \newpage
% \appendix

%%%%%%%%%%%%%%%
%%%%%%%%%%%%%%%
%%%%%%%%%%%%%%%
%%%%%%%%%%%%%%%

\section{Hardness of approximation.}\label{sec:hard}

In this {section} we prove Theorem~\ref{thm:approx},
which states that the sensitivity problem for (MBQP) is NP-hard to approximate.
Specifically, we show that computing an $(\alpha,\beta)$-approximation { of (\ref{eq:MBQPeq})} is NP-hard even if all variables are binary.
Our strategy is to create a trivial binary integer linear program which after changing one entry of $\bv$ by one captures a hard combinatorial property.
The hard combinatorial property we use are edge colorings of graphs.

Let $G = (V,E)$ be a simple graph. An edge coloring of $G$ is an assignment of colors to edges so that no incident edge will have the same color. The minimum number of colors required is called \textit{edge chromatic number} and denoted by $\chi'(G)$. The classical theory by Vizing~\cite{vizing1964estimate} states that:

\begin{theorem}[Vizing theorem] For any simple graph $G$, $\chi'(G) \in \{\Delta(G),\Delta(G)+1\}$ where $\Delta(G)$ is the maximum degree of vertices in $G$. 
\end{theorem}

Although Vizing theorem restricts edge chromatic number to two choices, it is still hard to distinguish between these two choices.
In fact, edge chromatic number is hard even in the special case of cubic graphs,
which are simple graph with every vertex having degree three.

\begin{theorem}[\cite{NP-Completeness}]
 It is NP-hard to determine the edge chromatic number of cubic graphs.
\end{theorem}

We can express edge coloring as an MBQP problem. Given a graph $G=(V,E)$ and an upper bound $H$ of its edge chromatic number $\chi'(G)$, then the classic formulation is
\begin{align*}
   & \min \sum_{i \in [H]} w_i \\
   & \text{s.t. } \sum_{i \in [H]} x_{ei} = 1,\forall e \in E \\
   & \ \ \ \ \ x_{r i} + x_{s i} \leq w_i,\forall r,s \in E, r \cap s \neq \emptyset , \forall i \in [H] \\
   & \ \ \ \ \ x \in \{0,1\}^{|E| \times [H]}, w \in \{0,1\}^{[H]}
\end{align*}
Here $w_i = 1$ means $i^{th}$ color is used and $x_{ri} = 1$ means edge $r$ is colored to be $i$. The first set of constraints requires that every edge must be colored by exactly one color. The second set of constraints requires that no adjacent edge will receive the same color.

When the given graph $G$ is a cubic graph, then $H = 4$ is an upper bound on $\chi'(G)$ by Vizing Theorem. Consider the following program:
\begin{align*}
   & z_1 : =\min \sum_{i \in [H]} w_i \\
   & \text{s.t. } \sum_{i \in [H]} - w_i \leq -4 \\
   & \ \ \ \ \ \sum_{i \in [H]} x_{ei} = 1,\forall e \in E \\
   & \ \ \ \ \ x_{r i} + x_{s i} \leq w_i,\forall r,s \in E, r \cap s \neq \emptyset , \forall i \in [H] \\
   & \ \ \ \ \ x \in \{0,1\}^{|E| \times [H]}, w \in \{0,1\}^{[H]}
\end{align*}

By Vizing Theorem, the optimal value $z_1$ is~$4$. When we change the first constraint from $ \sum_{i \in [H]} - w_i \leq -4$ to $ \sum_{i \in [H]} - w_i \leq -3$, we obtain the following program:
\begin{align*}
   & z_2 : =\min \sum_{i \in [H]} w_i \\
   & \text{s.t. } \sum_{i \in [H]} - w_i \leq -3 \\
   & \ \ \ \ \ \sum_{i \in [H]} x_{ei} = 1,\forall e \in E \\
   & \ \ \ \ \ x_{r i} + x_{s i} \leq w_i,\forall r,s \in E, r \cap s \neq \emptyset , \forall i \in [H] \\
   & \ \ \ \ \ x \in \{0,1\}^{|E| \times [H]}, w \in \{0,1\}^{[H]}
\end{align*}

{(\textbf{Proof of of Theorem~\ref{thm:approx}}):}
\begin{proof}
Let $z_1,z_2$ be the same as above. If the graph is $3$-edge-colorable, then $z_2 = 3$ and $|\Delta z| = 1$. Otherwise, $z_2 = 4$ and $|\Delta z| = 0$. In this case, any $({\Gamma_1},{\Gamma_2})$-approximation with ${\Gamma_2} \geq {\Gamma_1} > 0$ will return $p > 0$ if and only if the graph is $3$-edge-colorable. This implies that any such $({\Gamma_1},{\Gamma_2})$-approximation is NP-hard.
\end{proof}

%%%%%%%%%%%%%%%
%%%%%%%%%%%%%%%
%%%%%%%%%%%%%%%
%%%%%%%%%%%%%%%
%%%%%%%%%%%%%%%
\section{Duality theorem.}\label{sec:dual}
In this {section} we prove Theorem~\ref{thm:allstrong},
which characterizes when strong duality holds between (\ref{CP}) and (\ref{COP}). Assume that {a} lower bound {on} the optimal value (MBQP) is known in advance and { is denoted} by $l$. Recall that the feasible region of {the linear programming relaxation of} (MBQP) is 
$$\P = \{\xv :  \av_i^{\top} \xv  = b_i,\forall i \in [m], \xv \geq 0  \} \neq \emptyset.$$
In remaining part of this paper, we will introduce many new constants and functions. We present the following table to {reference} where these constants and functions are defined.

{
\begin{table}[h!]
\caption{{Notation used}}
% \label{tab:lpind}
\ra{0.5}
\setlength{\tabcolsep}{8pt}
\centering
 \begin{tabular}{cc} 
 \toprule
\text{name} & \text{source}
 \\  
  \midrule
$k$ & Lemma \ref{lemma_bounded_strict}  \\
$t_0$ & Lemma \ref{lem_lower_bound_t_bounded_case} \\
$t_1$ & Theorem \ref{thm_final_robust} \\
$t_2$ & Theorem \ref{thm_final_robust}  \\
$t_3$ & Proposition \ref{prop_threshold_no_binary_change} \\
$h(\cdot)$ & Remark \ref{rem_h_j} \\
$\mu(\cdot)$ & Remark \ref{rem_u_j} \\
$\eta$ & Theorem \ref{thm_cop_structure_unbound}\\
$\rho$ & Theorem \ref{thm_cop_structure_unbound} \\
 \hline
 \end{tabular}
\end{table}
}

This {section} is organized into several subsections.
Each of the parts of 
Theorem~\ref{thm:allstrong} is proved in a separate subsection.
{In the first subsection we present a proof of} 
Theorem~\ref{thm_final_robust},
which is used in the proof of Part~\ref{p2} of Theorem~\ref{thm:allstrong}.

%%%%%%%%%%%%%%%
%%%%%%%%%%%%%%%
\subsection{Theorem~\ref{thm:allstrong}\ref{p1}: Slater point when $\P$ is bounded.}\label{sec:boundedslater}

\begin{lemma}
\label{lemma_bounded_strict}
Let $H := T + \sum\limits_{i \in [m]} AA_i$. When $\P$ is bounded, $H$ is strictly copositive.  Therefore, there exists some number $k > 0$ that depends on $A$ such that $H - k I \in \COP$.
%and $H \underset{\COP}{\succeq} \coploweig I$ where $\coploweig > 0$ is a constant that depends on $A,\bv,\mathcal{B}$.
\end{lemma}

\begin{proof}
    $H$ is a sum of PSD {matrices} and therefore $H$ is clearly copositive.
    We prove that $H$ is strictly copositive by showing that for every non-zero $\yv := [t;\xv] \geq 0$, $\yv^{\top} H \yv > 0$. When $t > 0$, this is clearly true since $\yv^{\top} H \yv \geq \yv^{\top} T \yv = t > 0$. Thus, we may consider the case {where} $t = 0$. Suppose $H$ is not strictly copositive, then there exists some $\yv' := [0;\xv'] \geq 0$ ($\xv' \neq 0$) such that $(\yv')^{\top} H \yv' = 0$. This implies that
    \begin{align*}
    & (\yv')^{\top} AA_i \yv' = (\av_i^{\top} \xv')^2  = 0,\forall i \in [m] \implies \av_i^{\top} \xv' = 0, \forall i \in [m].
    \end{align*}
    {
    We first claim that $\xv'_j = 0,\forall j \in \mathcal{B}$. Suppose not and let $\xv^*$ satisfy $\xv \geq 0, \av_i^{\top} \xv = b_i,\forall i \in [m]$. It follows that $\xv^* + \pi \xv'$ still satisfies $\xv \geq 0, \av_i^{\top} \xv = b_i,\forall i \in [m]$ for any $\pi \geq 0$. However for sufficiently large $\pi$, there must exist some $j \in \mathcal{B}$ that either $(\xv^* + \pi \xv')_j < 0$ or $(\xv^* + \pi \xv')_j > 1$. This contradicts  Assumption \ref{eq:assumption}.} 
    
    {Now let $\xv''$ be a feasible solution of (\ref{eq:MBQPeq}). It follows that $\xv'' + \pi \xv'$ is a feasible solution of (\ref{eq:MBQPeq}) for any $\pi \geq 0$ since it satisfies all binary constraints and linear constraints. This implies that the feasible region of (\ref{eq:MBQPeq}) is unbounded and leads to a contradiction.
    }
    % \sout{Note that since by Assumption \ref{eq:assumption}, the constraints $\xv \in \P$ imply that $0 \leq x_j \leq 1$ for all $j \in \mathcal{B}$, we have that 
    % $\av_i^{\top} \xv' = 0, \forall i \in [m]$ implies $x_j' = 0,\forall j \in \mathcal{B}$.}
    % Thus, $\xv'$ is a non-zero recession direction of $\P$, implying that $\P$ is unbounded which leads to contradiction.
%    
%    To show $H \underset{\COP}{\succeq} \coploweig I$, consider the following quadratic program:
%    \begin{align*}
%        \coploweig := & \min_{\yv} \yv^{\top} H \yv \\
%        & \norm{\yv}_2 = 1 \\
%        & \yv \geq 0.
%    \end{align*}
%
%    Since $H$ is strictly copositive, $\coploweig > 0$. One can verify that $H \underset{\COP}{\succeq} \coploweig I$ by the definition of copositive matrix and $k$.
 \end{proof}

{\textbf{(Proof of Part~\ref{p1} of Theorem~\ref{thm:allstrong})}}:
\begin{proof}
    By Lemma \ref{lemma_bounded_strict}, $H$ is strictly copositive and therefore $C + \lambda H$ is strictly copositive for some sufficiently large $\lambda > 0$. By Slater condition, strong duality holds between (\ref{CP}) and (\ref{COP}). 
\end{proof}
%%%%%%%%%%%%%%%
%%%%%%%%%%%%%%%
\subsection{Theorem~\ref{thm_final_robust}: Local stability.}\label{sec:localstable}
%We proceed with the proof of (result.ii). The strategy is very similar to those proof of Lagrangian methods. We first verify that (MBQP) is very robust against tiny perturbation when $\P$ is bounded or $Q$ is PSD. Then  we construct several copositive matrices that behave like "Lagrangian penalties" of original constraints. Those results lead to the proof of (result.ii). It is worth mentioning that our proof is very natural and it explains why (CP-primal) is an exact reformulation of (MBQP) from the dual perspective.

%Consider (MBQP($\epsilon$)), 
%
%\begin{align}\tag{MBQP($\epsilon$)}
%    \begin{array}{rl}
%z(\epsilon) := & \min\limits_{\xv,\varepsilon^{(i)}}  \xv^{\top} Q \xv +  2\cv^{\top} \xv \\
%    &  \ \ \ \ \ \av_i^{\top} \xv  = b_i + \varepsilon^{(1)}_i,\forall i \in [m] \\
%    & \ \ \ \ \  x_j + \varepsilon^{(2)}_j \in \{0,1\}, \forall j \in \mathcal{B} \\
%    & \ \ \ \ \  \norm{\varepsilon^{(r)}}_{\infty} \leq \epsilon, \forall r \in \{1,2\}\\
%    & \ \ \ \ \  \xv \geq 0.
%    \end{array}
%\end{align}
For a fixed $\bv$, we will refer the feasible region of (MBQP($\epsilon$)) by $\mathcal{S}(\epsilon)$. $\mathcal{S}(\epsilon)$ is defined by some linear constraints and set of constraints $x_j + \varepsilon^{(2)}_j \in \{0,1\}, \forall j \in \mathcal{B}$. Since there are only finitely many choices of $x_j + \varepsilon^{(2)}_j$, $\mathcal{S}(\epsilon)$ can be viewed as a union of finitely many polyhedral. For any $\wv \in \{0,1\}^{\mathcal{B}}$, we define

\begin{align*}
& \mathcal{S}(\epsilon,\wv) :=\left\{ (\xv,{\varepsilon}) \;\middle\vert\;
   \begin{array}{@{}l@{}}
       \av_i^{\top} \xv = b_i + \varepsilon^{(1)}_i,\forall i \in [m] \\
       x_j + \varepsilon^{(2)}_j = w_j, \forall j \in \mathcal{B} \\
       {\varepsilon = [\varepsilon^{(1)} ; \varepsilon^{(2)}] \in \R^{m} \times \R^{|\mathcal{B}|},}  \\
       \norm{\varepsilon^{(r)}}_{\infty} \leq \epsilon, \forall r \in \{1,2\} \\
       \xv \geq 0
   \end{array}
\right\}. \\
& \zeta(\bv, \epsilon,\wv) := \min_{\xv,{\varepsilon}}\{ \xv^{\top} Q \xv + 2 \cv^{\top} \xv : {(\xv,{\varepsilon})} \in \mathcal{S}(\epsilon,\wv) \}. \\
& N(\epsilon) := \{ \wv \in \{0,1\}^{\mathcal{B}} : \mathcal{S}(\epsilon,\wv) \neq \emptyset \}. \\
& \xbar{N(\epsilon)} := \{0,1\}^{\mathcal{B}} \setminus N(\epsilon).
\end{align*}

Under this definition, we have
\begin{align*}
    & \mathcal{S}(\epsilon) = \bigcup\limits_{\wv \in \{0,1\}^{\mathcal{B}}} \mathcal{S}(\epsilon,\wv) = \bigcup\limits_{\wv \in N(\epsilon)} \mathcal{S}(\epsilon,\wv). \\
    & \zeta(\bv, \epsilon) = \min\limits_{\wv \in \{0,1\}^{\mathcal{B}}} \zeta(\bv, \epsilon,\wv) = \min\limits_{\wv \in N(\epsilon)} \zeta(\bv, \epsilon,\wv).
\end{align*}

We will prove that when $Q$ is PSD (including $Q =0$) or 
$\P$ is bounded, $\zeta(\bv, \epsilon)$ can be lower bounded by some linear function on $\epsilon$ if $\bv$ is fixed and  $\epsilon$ is small.
%\begin{theorem} (local stable theorem)
%\label{thm_final_robust} When $Q \text{ is PSD}$ or $\P$ is bounded,
%    there exists some $t_1,t_2$ that only depends only $A,\bv,\cv,Q,\mathcal{B}$ such that if $0 \leq \epsilon < t_1$, then then optimal value of (MBQP($\epsilon$)) is at least $l - \epsilon t_2$.
%\end{theorem}

To prove Theorem \ref{thm_final_robust}, we will establish several other statements.
Our main idea is to reduce obtaining lower bound on $\zeta(\bv, \epsilon)$ to finding the lower bound of finitely many quadratic programming problems. In particular, we will use Vavasis's result on the characterization of optimal solutions of general quadratic programming in \cite{VAVASIS199073} and show that the lower bound on $\zeta(\bv, \epsilon)$ can be viewed as piece-wise quadratic function on $\epsilon$ when $\epsilon$ is sufficiently small.

%When $Q =0$, Theorem \ref{thm_final_robust} is the special case of sensitivity of right-hand side change of mixed binary integer programming which is well studied in \textcolor{red}{schrijer textbook}. In our case, we use Vavasis's result to handle larger class. 

\begin{proposition}
\label{prop_bounded_eps}
    For any fixed $\epsilon \geq 0$, if $\P$ is bounded, then $\mathcal{S}(\epsilon)$ is bounded. Moreover, $\mathcal{S}(\epsilon_1)\subseteq \mathcal{S}(\epsilon_2)$ for all $\epsilon_2 \geq \epsilon_1 \geq 0$. 
\end{proposition}
\begin{proof}
    For any $\epsilon_1 \leq \epsilon_2$, if $[\xv;\varepsilon^{(i)}] \in \mathcal{S}({\epsilon_1})$, then clearly  $[\xv;\varepsilon^{(i)}] \in \mathcal{S}({\epsilon_2})$. This proves $\mathcal{S}(\epsilon_1)\subseteq \mathcal{S}(\epsilon_2)$. Moreover, $\mathcal{S}(\epsilon)$ can be viewed as a union of {feasible regions of finitely many linear programs of which right-hand-sides depend on $\epsilon$. The boundedness of $\P = \mathcal{S}({0})$ implies that each feasible region of linear programs (defined by choosing $\epsilon = 0$) is bounded.
    Since the boundedness of feasible regions of linear programs does not depdent on right-hand-sides, it follows that each feasible region of linear programs (defined by $\epsilon$) is bounded. Therefore, $\mathcal{S}({\epsilon})$ is bounded.}
    % it follows that boundedness of $\P = \mathcal{S}({0})$ implies boundedness of $\mathcal{S}({\epsilon})$. 
%     we prove the contrapositive statement, i.e., if $\mathcal{S}({\epsilon})$ is unbounded, then $\P$ is unbounded.
%     Suppose $\mathcal{S}({\epsilon})$ is unbounded, then consider its standard relaxation 
%     \[
% \mathcal{S}^{\text{relax}}(\epsilon) :=\left\{ (\xv \ \varepsilon^{(i)}) \;\middle\vert\;
%    \begin{array}{@{}l@{}}
% \av_i^{\top} \xv = b_i + \varepsilon^{(1)}_i,\forall i \in [m] \\
%     x_j + \varepsilon^{(2)}_j \in [0,1], \forall j \in \mathcal{B} \\
%     \norm{\varepsilon^{(r)}}_{\infty} \leq \epsilon, \forall r \in \{1,2\} \\
%        \xv \geq 0
%    \end{array}
% \right\}.
% \]
%     Since $\mathcal{S}^{\text{relax}}(\epsilon)$ is the relaxation of $\mathcal{S}(\epsilon)$, $\mathcal{S}^{\text{relax}}(\epsilon)$ is also unbounded.
%     Note $\mathcal{S}^{\text{relax}}(\epsilon)$ is defined by some linear constraints, this means there exists non-zero $[\xv_0;\varepsilon^{(i)}_0]$ such that $\av_i^{\top} \xv_0  = (\varepsilon^{(1)}_0)_i,\forall i \in [m],  (\xv_0)_j + (\varepsilon^{(2)}_0)_j = 0, \forall j \in \mathcal{B}, (\varepsilon^{(r)}_0)_t = 0, \forall r \in \{1,2\},\xv_0 \geq 0$. Rewriting those conditions, it yields there exists some non-zero $[\xv_0 ]$ such that $\av_i^{\top} \xv_0  =  0,\forall i \in [m], (\xv_0)_j = 0, \forall j \in \mathcal{B},\xv_0 \geq 0$. This implies that $\P$ is unbounded.
\end{proof}

\begin{proposition}    \label{prop_threshold_no_binary_change}
    There exists some threshold $t_3 > 0$ that only depends on $A,\bv,\mathcal{B}$ such that if $0 \leq \epsilon < t_3$, then $N(\epsilon) = N(0)$.
\end{proposition}
\begin{proof}
Observe that $\mathcal{S}(0,\wv) \subseteq \mathcal{S}(\epsilon,\wv)$ for all $\wv \in \{0,1\}^{\mathcal{B}}$. This implies $N(0) \subseteq N(\epsilon)$. 
    It suffices to show when $\epsilon$ is sufficiently small, for all $\wv \in \xbar{N(0)}$, we have that $\wv \in \xbar{N(\epsilon)}$. 
    For any $\wv \in \xbar{N(0)}$, consider the following linear programming:
    \begin{align}
        \begin{array}{rl}\label{local_lp}
 t_{\wv} :=&   \min \varphi \\
        &  \text{s.t. } b_i - \varphi \leq \av_i^{\top} \xv  \leq b_i + \varphi,\forall i \in [m] \\
    & \ \ \ \ \ x_j  = w_j,\forall j \in \mathcal{B} \\
     & \ \ \ \ \ \xv \geq 0, \varphi \geq 0.
        \end{array}
    \end{align}
    This linear program is clearly feasible by choosing $\varphi$ to be sufficiently large. Moreover, $t_{\wv} > 0$. Otherwise, this will imply that $\wv \in N(0)$. Since this linear program is  bounded from below and feasible, its optimal $t_{\wv}$ exists and $t_{\wv} > 0$ ($t_{\wv}$ {cannot} arbitrarily go to $0$ due to the attainability of linear programming). If $\epsilon >0$ is sufficiently small such that
    \begin{align*}
    \max_{i \in [m]} (1+\norm{\av_i}_1) \epsilon < t_{\wv},
    \end{align*}
    then we claim that $\wv \in \overline{N(0)}$ implies $\xbar{N(\epsilon)}$. Suppose not, $\wv \in N(\epsilon)$ and then $S(\epsilon,\wv) \neq \emptyset$.
    That is, there exists some $\Bar{\xv},\Bar{\varepsilon}^{(i)}$ such that
    \begin{align*}
        & \av_i^{\top} \Bar{\xv} + \Bar{\varepsilon}^{(1)}_i = b_i,\forall i \in [m] \\
        & \Bar{x}_j + \Bar{\varepsilon}^{(2)}_j = w_j, \forall j \in \mathcal{B} \\
        & \norm{\Bar{\varepsilon}^{(r)}}_{\infty} \leq \epsilon, r \in \{1,2\}\\
         & \Bar{\xv} \geq 0
    \end{align*}
    %$t_{\xv^*}$ can be interpreted as minimum change of right hand side to ensure $\xv^*$ feasible.
    Then we can construct a feasible solution $\xv^*$ for (\ref{local_lp}) where $x^*_j := \begin{cases}
        \Bar{x}_j & \text{ if } j \notin \mathcal{B} \\
        \Bar{x}_j + \Bar{\varepsilon}_j^{(2)} & \text{ if } j \in \mathcal{B}.
    \end{cases}$.
    In this case, for any $i \in [m]$, it follows
    \begin{align*}
        |\av_i^{\top} \xv^* - b_i|  & = |\av_i^{\top} \Bar{\xv} + (\sum_{j \in \mathcal{B}} \Bar{\varepsilon}_j^{(2)} (\av_i)_j) - b_i| \\
        & = |-\Bar{\varepsilon}^{(1)}_i + \sum_{j \in \mathcal{B}}\Bar{\varepsilon}_j^{(2)} (\av_i)_j |\\
        & \leq (1+\norm{\av_i}_1) \epsilon \\
        & < t_{\wv} .
    \end{align*}
    This contradicts the fact that $t_{\wv}$ is the optimal value of (\ref{local_lp}).
    
    Therefore to ensure $N(\epsilon) = N(0)$ , one can choose $t_3 := \frac{\min\limits_{\wv \in \overline{N(0)}} t_{\wv}}{\max\limits_{i \in [m]} (1+\norm{\av_i}_1)} > 0$.
\end{proof}

%\begin{remark}
%\label{rem_sdp_quadratic_zero}
%    For any $A \in \SS_+^n$ and $\xv \in \R^{n}$, if $\xv^{\top} A \xv = 0$, then $A \xv =0$.
%\end{remark}

\begin{proposition}
\label{prop_value_exists}
    When $Q$ is PSD or $\P$ is bounded, then $\zeta(\bv, \epsilon,\wv)$ exists for all $\epsilon \geq 0$ and $\wv \in N(\epsilon)$.
\end{proposition}

\begin{proof}
If $\P$ is bounded, then $\mathcal{S}(\epsilon)$ is bounded by Proposition \ref{prop_bounded_eps}. Since $\mathcal{S}(\epsilon,\wv) \subseteq \mathcal{S}(\epsilon)$, we have that $\mathcal{S}(\epsilon,\wv)$ is bounded as well. Thus $\zeta(\bv, \epsilon,\wv)$ is the optimal value of minimizing a quadratic over a compact set and therefore $\zeta(\bv, \epsilon,\wv)$ exists.

If $Q$ is PSD, we will use Vavasis' characterization of the optimal solution of quadratic programming in \cite{VAVASIS199073}. We would like to point out that when $Q = 0$, there is a simpler argument. When $Q=0$, $\zeta(\bv,\epsilon,\wv)$ is the optimal value of some linear program whose right-hand-side is parameterized by $\epsilon$. When $\epsilon = 0$, $\zeta(\bv, 0,\wv)$ exists and therefore this linear program is both dual feasible and primal feasible. When $\epsilon > 0$, the corresponding linear program is primal feasible as $\mathcal{S}(0,\wv) \subseteq \mathcal{S}(\epsilon,\wv)$. Moreover, the feasible region of the dual linear program remains the same and therefore it is dual feasible. In this case, $\zeta(\bv, \epsilon,\wv)$ exists.

When $Q \neq 0$ and $Q$ is PSD, suppose $\zeta(\bv,\epsilon,\wv)$ does not exist. By the result of \cite{VAVASIS199073},$\zeta(\bv,\epsilon,\wv)$ must diverge to negative infinity and there exists $\Bar{\xv}$, $\Bar{\varepsilon}^{(1)}$, $\Bar{\varepsilon}^{(2)},\dv_0,\dv_1,\dv_2$ such that for all large enough $t$,
$[\Bar{\xv} + t \dv_0; \Bar{\varepsilon}^{(1)} + t \dv_1;\Bar{\varepsilon}^{(2)} + t \dv_2]$ is feasible and has a decreasing objective function.
Since $\mathcal{S}(\epsilon,\wv)$ is defined by some linear constraints, this implies that $\dv_1 = 0,\dv_2 = 0,\av_i^{\top} \dv_0 = 0 ,\forall i \in [m] \text{ and } (\dv_0)_j = 0,\forall j \in \mathcal{B} \text{ and } \dv_0
 \geq 0$. Its objective function takes form of $(\Bar{\xv} + t \dv_0)^{\top}Q(\Bar{\xv} + t \dv_0) + 2\cv^{\top} (\Bar{\xv} + t \dv_0) = (\Bar{\xv})^{\top} Q \Bar{\xv} + 2t \dv_0^{\top} Q \Bar{\xv} + t^2 \dv_0^{\top} Q \dv_0 + 2\cv^{\top} \Bar{\xv} + 2 t \cv^{\top} \dv_0$. Since for all large $t$, the objective function is decreasing with large enough $t$ and $Q$ is PSD, it must be that $\dv_0^{\top} Q \dv_0 = 0$ and therefore $Q \dv_0 = 0$, since $Q$ is PSD. %by Remark \ref{rem_sdp_quadratic_zero}. 
 This further implies that $\cv^{\top} \dv_0 < 0$ since the objective function is decreasing. In this case, there exists some $\dv_0$ such that $Q \dv_0 = 0$ and $\cv^{\top} \dv_0 < 0$ and $\av_i^{\top} \dv_0 = 0,\forall i \in [m] \text{ and } (\dv_0)_j = 0,\forall j \in \mathcal{B} \text{ and } \dv_0 \geq 0$. Pick any $\xv \in \mathcal{S}(0,\wv)$, one can verify that $\xv + t \dv_0$ is feasible for all $t \geq 0$ and its objective value goes to negative infinity as $t$ goes to infinity.
This shows the $\zeta(\bv,0,\wv)$ diverges to negative infinity, which leads to a contradiction.

\end{proof}

% {Consider any quadratic program (possibly non-convex) of form
% \begin{align*}
%     \text{qv}(\gv) :=  & \min \yv^{\top} H \yv + 2 \dv^{\top} \yv \\
%     & \textbf{s.t. } M \yv \leq \gv.
% \end{align*}
% Vavasis proves the following theorem
% }

% Consider any quadratic program (possibly non-convex) of form
% \begin{align*}
%     \text{qv}(\gv) :=  & \min \yv^{\top} H \yv + 2 \dv^{\top} \yv \\
%     & \textbf{s.t. } M \yv \leq \gv.
% \end{align*}

% Vavasis proves the following theorem \cite{VAVASIS199073}
% \begin{theorem}

% \end{theorem}

% \textcolor{red}{--------------------------------------------------------}

Now we are {almost} ready to present a proof for Theorem \ref{thm_final_robust}. {As we mentioned earlier, our proof is built on Vavasis' characterization of optimal solution of quadratic programs.}
Consider any quadratic program (possibly non-convex) of form
\begin{equation}
\label{eq:general_quadratic_optimization}
    \begin{aligned}    
           \text{qv}(\gv) :=  & \min \yv^{\top} H \yv + 2 \dv^{\top} \yv \\
    & \textbf{s.t. } M \yv \leq \gv. 
    \end{aligned}
\end{equation}
\begin{theorem}
    \label{thm_Vavasis' characterization}
    (Vavasis' characterization of the optimal solution of quadratic programs \cite{VAVASIS199073}) If $\text{qv}(\gv)$ exists, then (\ref{eq:general_quadratic_optimization}) is equivalent to 
    a certain convex quadratic program:
\begin{align}
\label{eq:vavasis}\tag{$QP(\Tilde{M})$}
\begin{array}{rl}
    & \min_{\yv}  \yv^{\top} H \yv + 2\dv^{\top} \yv \\
    & \text{s.t. } \Tilde{M} \yv = \Tilde{\gv},
 %   & \ \ \ \ \ \ \yv^{\top} H \yv > 0, \forall \yv \text{ for all} \Tilde{M} \yv = \Tilde{g}
 \end{array}
\end{align}
where (i) $\Tilde{M}$ is a $(k+l)\times n$ matrix and first $k$ rows corresponding to some $k$
inequalities of original inequalities $M \yv \leq \gv$ satisfied exactly and the last $l$ rows corresponding to some entries of $\yv$ are zero and (ii) $H$ is positive definite when restricted to the special affine subspace defined by $\Tilde{M} \yv = \Tilde{g}$.
Moreover the convex program admits a unique solution $\yv_{\Tilde{M}}$ and its optimal value $\text{qv}_{\Tilde{M}}(\gv)$
where $\yv_{\Tilde{M}}$ is a linear function of $\gv$ and $\text{qv}_{\Tilde{M}}(\gv)$ is a quadratic function of $\gv$.
\end{theorem}
Note that there are finitely many possible choices for $\Tilde{M}$. For each choice of $\Tilde{M}$, $\yv_{\Tilde{M}}$ only satisfies some of the original constraints, and it is not necessarily feasible. We say $\Tilde{M}$ is \textit{good} if $\yv_{\Tilde{M}}$ is feasible and we say $\Tilde{M}$ is \textit{almost good} if $\yv_{\Tilde{M}}$ is infeasible. 
We denote the set of all good $\Tilde{M}$ by $\mathcal{M}_1$ and  denote the set of all almost good $\Tilde{M}$ by $\mathcal{M}_2$. Due to the definition of $\Tilde{M}$, both $\mathcal{M}_1,\mathcal{M}_2$ are finite sets and $\mathcal{M}_1 \cup \mathcal{M}_2$ is independent of $\gv$. %(However, each of them may depend on $\gv$). 
Theorem \ref{thm_Vavasis' characterization} states that optimal solution of (\ref{eq:general_quadratic_optimization}) will be some $\yv_{\Tilde{M}}$ for some $\Tilde{M}$ that is good. 

% \textcolor{red}{--------------------------------------------------------}

{\textbf{(Proof of of Theorem~\ref{thm_final_robust}}):}
\begin{proof}
By Proposition \ref{prop_threshold_no_binary_change}, when $\epsilon < t_3$ for some $t_3 > 0$ only depends on $A,\bv,\mathcal{B}$ , $N(\epsilon) = N(0)$. Pick any $\wv \in N(\epsilon)$, we would like to study $\zeta(\bv,\epsilon,\wv)$ which is the value function of a quadratic program whose right-hand side is parameterized by $\epsilon$. By Proposition \ref{prop_value_exists}, $\zeta(\bv,\epsilon,\wv)$ exists. 

In our case, we can express our program in inequality form to apply Theorem \ref{thm_Vavasis' characterization}:
\begin{align*}
    \begin{array}{rl}
& \min\limits_{\xv,\varepsilon^{(i)}}  \xv^{\top} Q \xv +  2\cv^{\top} \xv \\
    &  \ \ \ \ \  b_i - \varepsilon^{(1)}_i \leq \av_i^{\top} \xv  \leq b_i + \varepsilon^{(1)}_i,\forall i \in [m] \\
    & \ \ \ \ \   w_j - \varepsilon^{(2)}_j \leq x_j \leq w_j + \varepsilon^{(2)}_j , \forall j \in \mathcal{B} \\
    & \ \ \ \ \  0 \leq \varepsilon^{(r)}_t \leq \epsilon, \forall r \in \{1,2\}\\
    & \ \ \ \ \  \xv \geq 0.
    \end{array}
\end{align*}
{It is straightforward to see that the feasible region of the above program is the same as $S(\epsilon,\wv)$ where $\varepsilon$ is now signed.}
The right-hand-side of above program depends on $\epsilon \text{ and } \wv$. Let us refer to $\mathcal{M}_1,\mathcal{M}_2,\yv_{\Tilde{M}},qv_{\Tilde{M}}(\gv)$ corresponding to some $\epsilon \text{ and } \wv$ as $\mathcal{M}_1(\epsilon,\wv),\mathcal{M}_2(\epsilon,\wv),\yv(\epsilon,\wv,\Tilde{M}),\text{qv}(\epsilon,\wv,\Tilde{M})$. As our right hand side is a linear function of $\epsilon$, $\text{qv}(\epsilon,\wv,\Tilde{M})$ is a quadratic function of $\epsilon$ once $\wv$ and $\Tilde{M}$ are fixed. 

We fix some $\wv \in N(\epsilon)$ and prove that there exists some $t_4 >0 $ depending only on $A,\bv,\cv,Q,\mathcal{B}$, such that if $\epsilon < t_4$, then $\mathcal{M}_1(\epsilon,\wv) \subseteq \mathcal{M}_1(0,\wv)$. 
As pointed out earlier $\mathcal{M}_1(\epsilon,\wv) \cup \mathcal{M}_2(\epsilon,\wv) $ is independent of right hand sides and therefore is independent of $\epsilon$. Thus it suffices to prove that $\mathcal{M}_2(0,\wv) \subseteq \mathcal{M}_2(\epsilon,\wv)$. For $\Tilde{M} \in \mathcal{M}_2(\epsilon,\wv)$, $\yv_{\Tilde{M}}(\epsilon,\wv)$ is a continuous function of $\epsilon$ and therefore $\av^{\top}_i \yv_{\Tilde{M}}(\epsilon,\wv)$ is a continuous function of $\epsilon$. Thus, for sufficiently small values of $\epsilon$ if $\tilde{M} \in \mathcal{M}_2(0,\wv)$, then $\tilde{M} \in \mathcal{M}_2(\epsilon,\wv)$.

Thus when $\epsilon < \min\{t_3,t_4\}$, we have:
\begin{align*}
    \zeta(\bv,\epsilon) & = \min_{\wv \in N(\epsilon)} \min_{\Tilde{M} \in \mathcal{M}_1(\epsilon,\wv)} \text{qv}(\epsilon,\wv,\Tilde{M}) \\ & = \min_{\wv \in N(0)} \min_{\Tilde{M} \in \mathcal{M}_1(\epsilon,\wv)}  \text{qv}(\epsilon,\wv,\Tilde{M})\\ & \geq \min_{\wv \in N(0)} \min_{\Tilde{M} \in \mathcal{M}_1(0,\wv)} \text{qv}(\epsilon,\wv,\Tilde{M}) := \mathcal{T}(\epsilon).
\end{align*}

$\mathcal{T}(\epsilon)$ is the lower bound of $\zeta(\bv,\epsilon)$ with $\mathcal{T}(0) = \zeta(\bv,0) \geq l$. {Since there are only finitely many choices of $\wv \text{ and }\Tilde{M}$, $\mathcal{T}(\epsilon)$ is a piecewise quadratic function on $\epsilon$ with finitely many pieces. When $\epsilon \in  [0,\min\{t_3,t_4\}]$, $\mathcal{T}(\epsilon)$ is Lipschitz and therefore can be lower bounded by some affine function.}
That is, there exists some $t_2 > 0$ that depends on $A,\bv,\cv,Q,\mathcal{B}$, $\mathcal{T}(\epsilon) \geq l - t_2 \epsilon$ for $\epsilon \in [0,\min\{t_3,t_4\}]$ and this proves that $\zeta(\bv,\epsilon) \geq l - t_2 \epsilon$ for all $\epsilon \in [0,\min\{t_3,t_4\}]$.
\end{proof}

\subsection{Theorem~\ref{thm:allstrong}\ref{p2}: Constructing near optimal solution for (COP-dual). } \label{sec:constructsol}

In this subsection we prove part~\ref{p2} of Theorem~\ref{thm:allstrong} utilizing Theorem~\ref{thm_final_robust}. 
We first consider the case where the feasible set is bounded and then the case where it is unbounded and $Q$ is PSD.
Before that, we present some preliminary lemmas.

Given some number $f, g, l, \tau$, we remind the reader the two building blocks:
\begin{align*}
    & \KK_i  := b_i^2 T  -b_i A_i + AA_i,\forall i \in [m], \\
    & G_j(f,g,r)  := f\left(\sum\limits_{i \in [m]} \KK_i\right) - g N_j + r T,\forall j \in \mathcal{B}. 
 \end{align*} 
and $$H := T + \sum_{i \in [m]} AA_i.$$

The closed-form solution that we construct for (\ref{COP})  is of the form:
\begin{eqnarray}%\label{eq:dualsol}
U(f_1, f_2, g, r, \tau) = C + f_1\Bigl(\sum_{i = 1}^m \KK_i \Bigr)+ \sum_{j\in \mathcal{B}} G_j (f_2, g, r) + \tau H - l T.
\end{eqnarray}
{
We specify values for the parameters $f_1, f_2, g, r, \tau$ such that the above matrix is copositive and has objective value of $l - \epsilon$. 
Next consider the question of showing that $U(f_1, f_2, g, r, \tau)$ is copositive. It is easy to verify that {$\KK_i \succeq 0,H \succeq 0$ and therefore $\KK_i,H$ are copositive since the positive semidefinite cone is the inner approximation of copositive cone}. We also show that for sufficiently large $f, g, r$ we have that $G_j(f, g, r) \in \COP$ {in Lemma \ref{lem_small_copostive block}}. However, due to presence of the terms $C$ and $-lT$ in $U(f_1, f_2, g, r, \tau)$, one has to additionally verify its copositivity. Consider a non-negative vector $\yv := [t; \xv] \geq 0$ and we need to verify that $\yv^{\top}U(f_1, f_2, g, r, \tau) \yv \geq 0$. In the case when $t =0$, this follows from the fact that $Q \succeq 0$ or $H$ is strictly copositive when $\P$ is bounded {(see the proof of Theorem \ref{thm_cop_structure} and Theorem \ref{thm_cop_structure_unbound} for details)}.
% The case when $t = 0$ follows from the fact that $Q, \KK_i \succeq 0$ and $G_j(f, g, r), H \in \COP$. 
In the case when $t >0$, the building blocks $\KK_i$'s and $G_j(f, g, r)$'s behave like \emph{augmented Lagrangian penalties}  (see~\cite{gu2020exact,feizollahi2017exact} for strong duality results for general mixed integer convex quadratic programs) of the original constraints of the MBQP, implying the required inequality.
% \st{
% , i.e., if $\frac{1}{t} \xv$ is not feasible for MBQP($\epsilon$) then $\yv^{\top}U(f_1, f_2, g, r, \tau) \yv$ becomes large positive value in the following fashion:}
{In particular, given $\epsilon > 0$, we partition $\yv$ in the following cases and assert that $\yv^{\top}U(f_1, f_2, g, r, \tau)\yv \geq 0$ {(see the proof of Theorem \ref{thm_cop_structure} and Theorem \ref{thm_cop_structure_unbound} for details)}:}
\begin{itemize}
\item If $|\av_i^T \xv\frac{1}{t} - b_i| > \epsilon$, then it is easy to see that  $y^{\top}\KK_iy  = |\av_i^Tx - t\bv_i|^2$ and this ``penalty" yields that $\yv^{\top}U(f_1, f_2, g, r, \tau)\yv \geq 0$ for the selected values of parameters $(f_1, f_2, g, r, \tau)$.
\item If $\frac{x_j}{t} \in (\epsilon, 1- \epsilon) \cup (1+ \epsilon, \infty)$ (i.e., $\frac{x_j}{t}$ is far from being binary), then $y^{\top}G_j(f,g,r)y \gg 0$ and this ``penalty" yields that $\yv^{\top}U(f_1, f_2, g, r, \tau)\yv \geq 0$ for the selected values of parameters $(f_1, f_2, g, r, \tau)$.
\item {The remaining case is when $\frac{1}{t}\xv$ is ``almost" feasible for MBQP, where we refer such region as MBQP($\epsilon$). This case is handled by Theorem \ref{thm_final_robust}}
\end{itemize}
}

Also let:
 \begin{align*}  
    & U_0 : = U_0(\tau,l) := C  + \left(\sum_{i=1}^{m}  \frac{|\lambda_{\min}(Q)|+1}{k} \KK_i \right) +  (\tau H) +  (-lT),
\end{align*}
where $k$ is defined in Lemma~\ref{lemma_bounded_strict}.
%All of $H, \KK_i,G_j(f,g,r),U_0(\alpha,l)$ will be building blocks of our closed-form construction. As we mentioned earlier, we will give a Lagrange proof.
As mentioned earlier, $\KK_i$ will serve as a penalty block of $\av^{\top}_i \xv = b_i$ and $G_j(f,g,r)$ will serve as a penalty block of $x_j \in \{0,1\}$. Consider $\yv := [t;\xv] \in \R^{n+1}_+$, it follows
\begin{align*}
    & {\yv^{\top}( \KK_i) \yv}  = (b_i t - \av_i^{\top} \xv )^2,\forall i \in [m], \\
   & {\yv^{\top}(- N_j) \yv} = 2   x_j (t - x_j),\forall j \in \mathcal{B}.
\end{align*}
Note that all of $H, \KK_i$ are PSD matrices and therefore copositive. To see that $\KK_i$ is PSD, note it can be written as
    \begin{align*}
        \KK_i = \begin{bmatrix}
            b_i^2 & -b_i (\av)_i^{\top} \\ -b_i \av_i & \av_i (\av)_i^{\top} 
        \end{bmatrix} = \left(\begin{bmatrix}
            -b_i \\ \av_i
        \end{bmatrix}\right) \left(\begin{bmatrix}
            -b_i \\ \av_i
        \end{bmatrix}\right )^{\top}.
    \end{align*}
The matrix $G_j(f,g,r)$ is parameterized by $f,g,r$. It is {not} necessarily copositive for arbitrary choice of $f,g,r$ and $G_j(f,g,r)$ has objective value equals to $r$. The next lemma proves that for any positive $r$ and positive $g$, one can choose {sufficiently} large $f$ so that $G_j(f,g,r)$ is copositive.

In the remaining of this {{section}}, for any matrix $A \in \SS^{n+1}$, we will use $(A)_{[\xv]}$ to refer to principal submatrix of $A$ corresponding to the indices of $[\xv]$. We will refer other principal submatrices like $(A)_{[x_j]}$ in the similar manner.
Recall that we make the following assumption (\ref{eq:assumption}):
\begin{equation*}
\xv \geq 0,\  \av_i^{\top} \xv = b_i,\ \forall i \in[m]
\quad\implies\quad
0 \leq x_j \leq 1 \text{ for all } j \in \mathcal{B}.
\end{equation*}
We first establish several remarks related to the assumption (\ref{eq:assumption}). 
\begin{remark}
\label{rem_h_j}
    Fix some index $j \in \mathcal{B}$ and some positive number $\eta > 0$. Then for any $\xv \geq 0 \text{ such that } x_j \geq \eta $, there exists some strict positive number $h_j(\eta) > 0$ that only depends on $A,b,\eta,j$ such that 
    $$\max\limits_{i \in [m]} |\av_{i}^{\top} \xv| \geq h_j(\eta).$$
\end{remark}

\begin{proof}
    Consider the following linear programming:
    \begin{align}
    \label{LP_in_rem_h_j}
    \begin{array}{rl}
        h(\eta) := & \min h \\
        & \text{s.t.}  -h \leq \av^{\top}_i \xv \leq h,\forall i \in [m], \\
        & \ \ \ \ \ \xv \geq 0, h \geq 0, \\
        & \ \ \ \ \ x_j \geq \eta.
    \end{array}
    \end{align}
    (\ref{LP_in_rem_h_j}) is clearly feasible by choosing sufficiently large $h$. Moreover, $h(\eta) > 0$. Otherwise, if $h(\eta)  = 0$, let $\xv^*,h^*$ be the corresponding optimal solution. Note that $\xv^* \neq 0$. Pick $\xv_0 \in \P$, it follows that $\xv_0 + \lambda \xv^* \in \P,\forall \lambda \geq 0$. With sufficiently large $\lambda$, the assumption (\ref{eq:assumption}) for $j \in \mathcal{B}$ will be violated. 
\end{proof}

Similarly, one can prove the following remarks using the similar approach :
\begin{remark}
\label{rem_u_j}
    Fixed $j \in \mathcal{B}$ and $\eta > 0$. Then for any $\xv \geq 0 \text{ such that } x_j \geq 1 + \eta $, there exists some strict positive number $u_j(\eta) > 0$ that only depends on $A,b,\eta,j$ such that 
    $$\max\limits_{i \in [m]} |\av_{i}^{\top} \xv - b_i| \geq u_j(\eta).$$
\end{remark}

\begin{remark}
    \label{rem_key_assumpton}
    For any $\xv \geq 0$, if $\av_i^{\top} \xv = 0,\forall i \in [m]$, then $x_j = 0,\forall i \in \mathcal{B}$.
\end{remark}

\begin{lemma}
    \label{lem_small_copostive block}
    For any $j \in \mathcal{B},g > 0, r > 0$, let $h_j(\eta)$ be as defined in Remark \ref{rem_h_j}. There exists some strictly positive number ${p}_j > 0$ that only depends on $A,\bv,r,g,j$ such that
    if $f$ satisfies
    \begin{eqnarray}
    \label{eq:choice_of_f_small_block}
    f \geq \max \left\{ \frac{2g}{h_j(1)^2}, \frac{2g}{{p_j}^2}, \frac{g^2 + 2rg}{r {p_j}^2} \right\},
    \end{eqnarray}

    then $G_j(f,g,r) = (f\sum\limits_{i \in [m]} \KK_i) - g N_j + r T$ is copositive. 
\end{lemma}
\begin{proof}
    For any $ \yv = [t;\xv] \geq 0$, we consider several cases of $\yv$ and assert that 
    $$\Lambda := \yv^{\top} ((f\sum\limits_{i \in [m]} \KK_i) - g N_j + r T) \yv \geq 0.$$ 
    Note that since both $\sum\limits_{i \in [m]} \KK_i$ and $T$ are PSD, the only term that could potentially make $\Lambda$ negative is $-g N_j$. 
    
    If $t = 0$, note that $(- g N_j)_{[\xv]}$ is a diagonal matrix with one negative entry such that $(- g N_j)_{[x_j]} = -2g$ and other entries all zeros. Therefore, if $x_j = 0$, then $\Lambda \geq 0$. If $x_j \neq 0$, we can further assume $x_j = 1$ after scaling properly. Let $h_j(1) > 0$ be the constant in Remark \ref{rem_h_j}. Then it follows that
        \begin{align*}
        \Lambda & = \yv^{\top} \left((f\sum\limits_{i \in [m]} \KK_i) - g N_j + r T\right) \yv \\
                & = \xv^{\top} \left(f \sum_{i \in [m]} \av_i \av_i^{\top} \right) \xv - 2g \\
                & = f \left(\sum_{i \in [m]} (\av_i^{\top} \xv)^2 \right) - 2 g \\
                & \geq f h_j(1)^2 - 2g & (\text{by Remark } \ref{rem_h_j}) \\
                & \geq 0. & (\text{by } (\ref{eq:choice_of_f_small_block}))
    \end{align*}
    If $t \neq 0$, we may assume $t=1$ after scaling properly. If  $t = 1 \geq x_j$, ${\yv}^{\top}(- g N_j) {\yv}^{\top} 
 \geq 0$. Since all $\KK_i$ and  $T$ are PSD, then $\Lambda \geq 0.$ Thus it remains to consider the case when $t < x_j,t = 1$ and let $v := x_j -t = x_j - 1 > 0$. If $v \leq \min\{\frac{r}{4g},1\}$, then it follows that
 \begin{align*}
        \Lambda & = \yv^{\top} \left((f\sum\limits_{i \in [m]} \KK_i) - g N_j + r T\right) \yv \\
                & \geq \yv^{\top} \left(- g N_j + r T \right) \yv \\
                & = -2 g (v+1)v + r \\
                & \geq -4 g v + r & (0 \leq v \leq 1) \\
                & \geq 0. &(v\leq \frac{r}{4g})
    \end{align*}
    So we may assume that $v \geq \min\{\frac{r}{4g},1\}$. In this case, we show that at least one of $|\av^{\top} \xv_i - b_i|$ is considerably large and increases at least linearly with respect to $v$. To be more precise, we are considering the following linear fractional programming:
\begin{align*}
  % \label{pj}
  \begin{array}{rl}
        p_j := & \min\limits_{\xv,\psi,v} \frac{\psi}{v} \\
        & \text{s.t. } b_i - \psi \leq \av_i^{\top} \xv \leq b_i + \psi,\forall i \in [m], \\
        & \ \ \ \ \ x_j \geq 1 + v, \\
        & \ \ \ \ \ v \geq \min\{\frac{r}{4g},1\}, \\
        & \ \ \ \ \ \xv \geq 0, \psi \geq 0.
  \end{array}
\end{align*}
    This linear fractional programming can be exactly formulated as the following linear programming, which is known as Charnes-Cooper transformation \cite{Cooper1962}.
\begin{align}
  \label{pj_tilda}
  \begin{array}{rl}
        {p_j} := & \min\limits_{\Tilde{\xv},\Tilde{\psi},\Tilde{v},s} \Tilde{\psi} \\
        & \text{s.t. } b_i s - \Tilde{\psi} \leq \av_i^{\top} \Tilde{\xv} \leq b_i s + \Tilde{\psi},\forall i \in [m], \\
        & \ \ \ \ \ \Tilde{\xv}_j \geq s + \Tilde{v}, \\
        & \ \ \ \ \ \Tilde{v} \geq \min\{\frac{r}{4g},1\} s, \\
        & \ \ \ \ \ \Tilde{\xv} \geq 0, s \geq 0, \Tilde{\psi} \geq 0, \Tilde{v} = 1.
  \end{array}
\end{align}
% One can prove that $p_j = \Tilde{p_j}$ and this is known as Charnes-Cooper transformation \cite{Cooper1962}.
% We can see that any feasible solution $(\xv,\psi,v)$ in $(\ref{pj})$ can be mapped to a feasible solution in $(\ref{pj_tilda})$ by setting $\Tilde{\xv}_j = \frac{\xv}{v},\Tilde{\psi} = \frac{\psi}{v},\Tilde{v} = \frac{v}{v},s = \frac{1}{v}$ without changing objective value and this implies that $p_j \geq \Tilde{p_j}$ (one can actually further show that $p_j = \Tilde{p_j}$). 
We can further claim that ${p_j} > 0$ otherwise if ${p_j} = 0$ and let $\Tilde{\xv^*},\Tilde{\psi^*},\Tilde{v^*},s^*$ be the corresponding optimal solution of (\ref{pj_tilda}). If $s^* = 0$, then $\av_i^{\top} \Tilde{\xv^*} = 0,\forall i \in [m]$. By Remark \ref{rem_key_assumpton}, this implies that $\Tilde{\xv_j} = 0$. However, this contradicts that $\Tilde{\xv}_j \geq s + \Tilde{v} \geq 1$. In the case when $s^* > 0$, one can again verify that $\frac{\Tilde{\xv^*}}{s^*}$ will violate  assumption (\ref{eq:assumption}) by using Remark~\ref{rem_u_j}. Finally, it follows that
    \begin{align*}
        \Lambda & = \yv^{\top} \left(\left(f\sum\limits_{i \in [m]} \KK_i\right) - g N_j + r T\right) \yv \\
        & = f \left(\sum_{i\in[m]} (\av_i^{\top} \xv - b_i)^2 \right) - 2g(v+1)v + r \\
        & \geq f {p_j}^2 v^2 - 2g(v+1)v + r & (\text{by the definition of } {p_j})  \\
        & = (f {p_j}^2 - 2g) v^2 - 2gv + r \\
        & \geq 0. & (\text{by } (\ref{eq:choice_of_f_small_block}))
    \end{align*}
The last inequality comes from the fact that the lower bound $(f {p_j}^2 - 2g) v^2 - 2gv + r)$ is a quadratic function with respect to $v$. If $f {p_j}^2 - 2g \geq 0$ and $g^2 - r(f {p_j}^2 - 2g)\leq0$, then this lower bound is non-negative for all $v$ and this condition is guaranteed by our choice of $f$.
\end{proof}

\begin{remark}
\label{rem_ensure_strict_copostive}
    Consider two symmetric matrices $A,B$ such that $B \underset{\COP}{\succeq} \eta I$ for some $\eta > 0$, then $A + \frac{|\lambda_{\min}(A)|+1}{\eta} B$ is strictly copositive.
\end{remark}
\begin{proof}
    {By the definition of $\lambda_{\min}$, we have $\yv^{\top} A \yv \geq \lambda_{\min}(A)$ for all $\yv$ such that $\norm{\yv}_2 = 1$. Since $B \underset{\COP}{\succeq} \eta I$, this implies that $\yv^{\top} B \yv \geq \eta$ for all $\yv \geq 0$ such that $\norm{\yv}_2 =1$. In this case, for all $\yv \geq 0$ such that $\norm{\yv}_2 = 1$, it follows that
    \begin{align*}
        \yv^{\top} (A + \frac{|\lambda_{\min}(A)|+1}{\eta} B) \yv & =  \yv^{\top} A \yv + \yv^{\top} (\frac{|\lambda_{\min}(A)|+1}{\eta} B) \yv \\
        & \geq \lambda_{\min}(A) + |\lambda_{\min}(A)|+1 \\
        & > 0.
    \end{align*}
    This therefore implies that $A + \frac{|\lambda_{\min}(A)|+1}{\eta} B$ is strictly copositive.}
\end{proof}

\begin{lemma}
\label{lem_lower_bound_t_bounded_case} Recall $U_0 = U_0(\tau,l) = C  + \left(\sum_{i=1}^{m}  \frac{|\lambda_{\min}(Q)|+1}{k} \KK_i\right) +  (\tau H) +  (-lT)$.
    If $\P$ is bounded, then there exists some $t_0 > 0$ that only depends on $A,\bv,\cv, Q, l, \tau$
    such that for any $\yv := [t ; \xv] \geq 0$ with $t < t_0$ and $\norm{[t ; \xv]}_2 = 1$, $\yv^{\top} U_0 \yv \geq 0$.
\end{lemma}

\begin{proof}
    We first show that $(U_0)_{[\xv]}$ is strictly copositive. By Lemma \ref{lemma_bounded_strict}, we have $H \underset{\COP}{\succeq} \coploweig I$. This implies that $(H)_{[\xv]} \underset{\COP}{\succeq} \coploweig I$.
    As $(H)_{[\xv]} = \sum\limits_{i \in [m]} \av_i(\av_i)^{\top} $,
this means that $\sum\limits_{i \in [m]} \av_i(\av_i)^{\top} \underset{\COP}{\succeq} \coploweig I$. Applying Remark $\ref{rem_ensure_strict_copostive}$, we obtain that $(U_0)_{[\xv]} = Q + \frac{|\lambda_{\min}(Q)|+1}{k} (\sum\limits_{i \in [m]} (\av_i\av_i^{\top})) + \tau H_{[x]} $ is strictly copositive. 
Consider the following program:
    \begin{align*}
        t_0 := & \min t \\
        & \text{s.t. }  \yv  = [t ; \xv] \geq 0, \\
& \ \ \ \ \ \norm{\yv}_2 = 1, \\
& \ \ \ \ \ \yv^{\top} U_0 \yv  \leq 0.
    \end{align*}
    If this program is infeasible, we may set $t_0 = \infty$.
    Since the feasible region is a compact set and if this program is  feasible, the optimal value exists and is attained. 
    We first claim that $t_0 > 0$. Otherwise,
    there exits some $\yv' := [0;\xv'] \geq 0$ such that $(\yv')^{\top} U_0 \yv' \leq 0$. This contradicts the fact that $(U_0)_{[\xv]}$ is strictly copositive. In this case, 
   the definition of $t_0$ guarantees that 
   % $$\forall \yv := [t ; \xv] \geq 0, t < t_0, \norm{[t ; \xv]}_2 = 1 \implies  \yv^{\top} U_0 \yv > 0.$$
   for any $\yv := [t;\xv] \geq 0 $ with $t < t_0$ and $\norm{[t ; \xv]}_2 = 1$, $\yv^{\top} U_0 \yv > 0$. 
\end{proof}

We proceed to prove Part~\ref{p2} of Theorem~\ref{thm:allstrong} in the case that the feasible set is bounded.

\begin{theorem} \label{thm_cop_structure} (Part~\ref{p2} of Theorem~\ref{thm:allstrong}- bounded case) 
Let $l$ be the lower bound of the optimal value of (MBQP) that is $l \leq z(\bv)$ and $l_+ := \max\{l,0\}$. Let $k$ be the constant defined in Lemma \ref{lemma_bounded_strict}, let $t_1,t_2$ be as in Theorem \ref{thm_final_robust} and let $u_j(\cdot)$ be the constant defined in Remark \ref{rem_u_j}. When $\P$ is bounded, for any $\epsilon_0 \in (0,t_1)$ and any $r > 0$,
\begin{align*}
    U := U(f_1, f_2, g, r, \tau, l) := C + f_1 \left(\sum_{i \in [m]} \KK_i \right) + \sum_{j \in \mathcal{B}} G_j(f_2, g, r) + \tau H - l T
\end{align*}
is copositive for $f_1,f_2,g,\alpha$ satisfying the following rules:
\begin{enumerate}[label=(\text{rule}.\roman*),align=left]
    \item $f_2, g, r$ satisfies the condition in Lemma \label{rule:1} \ref{lem_small_copostive block} so that $G_j(f_2, g, r)$ is copositive,
    \item $\tau = \frac{\epsilon_0 t_2}{k}$, \label{rule:2}
    \item $t_0$ is defined in Lemma \ref{lem_lower_bound_t_bounded_case}, which depends on $\tau, l$ and therefore depends on $\epsilon_0,l$, \label{rule:3} 
    \item $f_1 \geq \max \left\{ \frac{|\lambda_{\min}(Q)|+1}{k}, \frac{-\lambda_{\min}(C) + l_+}{t_0^2 \epsilon^2_0}, \max\limits_{j \in \mathcal{B}}\{ \frac{\lambda_{\min}(C) + l_+}{t_0^2 u_j(\epsilon_0)^2} \}\right\}$, \label{rule:4}
    \item $g \geq \frac{\lambda_{\min}(C) - l_+}{2\epsilon^2_0t_0^2}$. \label{rule:5}
\end{enumerate}
Moreover, {for the given $\epsilon_0,r,l$, the objective value of $U$ is $l- t_2 \epsilon_0 \cdot (1+\sum\limits_{i \in [m]} b_i^2) - r\cdot|\mathcal{B}|$.} This objective value can be arbitrarily {close} to $l$ as $r$ and $\epsilon_0$ goes to zero.
\end{theorem}

We would like to point out that our bound in Theorem \ref{thm_cop_structure} is rather loose, and we are only seeking sufficient conditions to ensure $U$ is copositive. 
We begin to prove Theorem \ref{thm_cop_structure}.
\begin{proof}

By \ref{rule:4}, it follows that
\begin{align*}
    f_1 - \frac{|\lambda_{\min}(Q)|+1}{k} \geq 0.
\end{align*}
This further implies
\begin{align*}
    U - U_0 = \left(f_1 - \frac{|\lambda_{\min}(Q)|+1}{k}\right) \left(\sum_{i\in [m]} \KK_i\right) + \sum_{j \in \mathcal{B}} G_j(f_2, g, r) \underset{\COP}{\succeq} 0.
\end{align*}
Combining with Lemma \ref{lem_lower_bound_t_bounded_case}, this implies that there exists some $t_0 > 0$ such that any $\yv := [t ; \xv] \geq 0$ with $t < t_0$ and $\norm{[t ; \xv]}_2 = 1$, $\yv^{\top} U \yv \geq 0$.
Thus we can assume that $t \geq t_0$ and define the following two sets:
{
\begin{align}
    & \I(\epsilon_0)  := \left\{ (t,\xv) \geq 0 :  \frac{x_j}{t} \in [0,\epsilon_0] \cup [1-\epsilon_0,1+\epsilon_0],\forall j \in \mathcal{B} \right\}, \label{eq:idefn}\\
    & \P(\epsilon_0)  := \left\{ (t,\xv) \geq 0 : \frac{1}{t} (\av_i^{\top} \xv) \in [b_i-\epsilon_0,b_i+\epsilon_0],\forall i \in [m]  \right\}. \label{eq:pdefn}
\end{align}
}

One interpretation of such sets is that $\I(\epsilon_0)$ is the approximate version of $0$-$1$ integrality and $\P(\epsilon_0)$ is the approximate version of $\P$. 

Now assume $\yv \geq 0 := [t;\xv]$ with $\norm{\yv}_2 = 1$ and $t \geq t_0$. We are going to prove that for all such $\yv \geq 0,\yv^{\top} U \yv \geq 0$.
We consider three cases:
\begin{enumerate}[label=(\text{case}.\arabic*),align=left]
    \item $[t ; \xv] \in \xbar{\P(\epsilon_0)}$
    \item $[t ; \xv] \in \P(\epsilon_0) \cap \I(\epsilon_0)$
    \item $[t ; \xv] \in \P(\epsilon_0) \cap \xbar{\I(\epsilon_0)}$
\end{enumerate}
Recall that
{
\begin{align*}
   \yv^{\top} U \yv & = \yv^{\top} C \yv + \sum_{i \in [m]} 
 \underbrace{\yv^{\top} f_1 (\KK_i) \yv}_{ \geq 0} + \sum_{j \in \mathcal{B}} \underbrace{\yv^{\top} ( G_j(f_2,g,r))\yv}_{\geq 0} + \underbrace{\yv^{\top}\tau H \yv}_{\geq 0} - \yv^{\top} l T \yv
\end{align*}
}
The only term that only potentially makes $\yv^{\top} U \yv$ negative is $\yv^{\top} C \yv + \yv^{\top} (-lT) \yv$. Given $\norm{\yv}_2 = 1$, a (trivial) lower bound on $\yv^{\top} C \yv + \yv^{\top} (-lT) \yv$ is
\begin{align}
\label{trivial_lower_bound}
    \yv^{\top} C \yv + \yv^{\top} (-lT) \yv & = \yv^{\top} C \yv - t^2 l 
     \geq \lambda_{\min}(C) - l_+.
\end{align}
\text{(case.1)} 
If $[t;\xv] \in \xbar{\P(\epsilon_0)}$, there must exist an index $e \in [m]$  such that
\begin{align}
    \label{eq:case1_bounded}
    \left|\frac{1}{t} (\av_e^{\top} \xv)  - b_e \right| > \epsilon_0.
\end{align}
It follows that
\begin{align*}
    \yv^{\top} U \yv & \geq \yv^{\top} C \yv + \yv^{\top} (-lT) \yv + \yv^{\top} \KK_i \yv \\
    & = \yv^{\top} C \yv + \yv^{\top} (-lT) \yv + t^2 f_1 \left(\frac{1}{t}(\av_e^{\top} \xv + s_e) - b_e\right)^2 \\
    & \geq \lambda_{\min}(C) - l_+ + t_0^2 f_1 \epsilon^2_0 \ & (\text{by } (\ref{trivial_lower_bound}),(\ref{eq:case1_bounded}))   \\
    & \geq 0 & (\text{by } \ref{rule:4}).
\end{align*}
\text{(case.2)}  
If $[t;\xv] \in \P(\epsilon_0) \cap \I(\epsilon_0)$, then $\frac{1}{t}(\xv)$ is a feasible solution of MBIP($\epsilon_0$). Since $\epsilon_0 < t_1$ by our choice, Theorem \ref{thm_final_robust} implies that
\begin{align}
\label{eq:case2_bounded}
   \left(\frac{1}{t}\xv\right)^{\top} Q \left(\frac{1}{t}\xv \right) +  2\cv^{\top} \frac{1}{t}\xv   \geq l - t_2 \epsilon_0. 
\end{align}
Therefore, it follows that
\begin{align*}
\yv^{\top} U \yv & \geq \yv^{\top} C \yv + \yv^{\top} (-lT) \yv +  \yv^{\top} \tau H \yv \\
& = \xv^{\top} Q \xv + 2 t \cv \xv + t^2 (-l) + \yv^{\top} \tau H \yv \\
& = t^2 \left(\left(\frac{1}{t}\xv\right)^{\top} Q \left(\frac{1}{t}\xv\right) +  2\cv^{\top} \frac{1}{t}\xv-l \right) + \yv^{\top} \tau H \yv  \\
& \geq t^2 (l - t_2 \epsilon_0  - l) + \yv^{\top} \tau H \yv & (\text{by } (\ref{eq:case2_bounded})) \\
& \geq (-t_2 \epsilon_0) + \tau k & (H \text{ is strictly copositive and } t \leq 1 ) \\
& = 0 & (\text{by } \ref{rule:3}).
\end{align*}
\text{(case.3)} 
If $[t;\xv] \in \P(\epsilon_0) \cap \xbar{\I(\epsilon_0)}$, this implies that there exists $e \in \mathcal{B}$ such that
$$ \text{either } \frac{x_e}{t} \in (1+\epsilon_0,\infty) \text{ or } \frac{x_e}{t} \in (\epsilon_0,1-\epsilon_0).$$
If $\frac{x_e}{t} \in (1+\epsilon_0,\infty)$, by Remark \ref{rem_u_j}, there exists some $w \in [m]$, some $e \in \mathcal{B}$ and some $u_e(\epsilon_0) > 0$  such that 
\begin{align}
\label{eq:case3_bounded}
    \left|\frac{1}{t} (\av_w^{\top} \xv)  - b_w \right| > u_e(\epsilon_0)
\end{align}
and it follows that
\begin{align*}
    \yv^{\top} U \yv & \geq \yv^{\top} C \yv + \yv^{\top} (-lT) \yv + \yv^{\top} \KK_w \yv \\
    & = \yv^{\top} C \yv + \yv^{\top} (-lT) \yv + t^2 f_1 \left(\frac{1}{t}(\av_w^{\top} \xv ) - b_w \right)^2 \\
    & \geq \lambda_{\min}(C) - l_+ + t_0^2 f_1  {u_e(\epsilon_0)}^2 \ & (\text{by } (\ref{trivial_lower_bound}),(\ref{eq:case3_bounded}))   \\
    & \geq 0 & (\text{by } \ref{rule:4}).
\end{align*}
If $\frac{x_e}{t} \in (\epsilon_0,1-\epsilon_0)$, this implies that
\begin{align}
    \label{eq:case3_bounded_last}
    x_e \geq \epsilon_0 t \geq \epsilon_0 t_0 \text{ and } t - x_e  \geq \epsilon_0 t \geq \epsilon_0 t_0.
\end{align}
% \begin{align*}
%     & x_e \geq \epsilon_0 t \geq \epsilon_0 t_0 \\
%     & t - x_e  \geq \epsilon_0 t \geq \epsilon_0 t_0
% \end{align*}
Therefore, it follows that
\begin{align*}
   \yv^{\top} U \yv  & \geq \yv^{\top} C \yv + \yv^{\top} (-lT) \yv + \yv^{\top} ( G_e(f_2,g,r))\yv \\
   & \geq \lambda_{\min}(C) - l_+ + \yv^{\top} ( G_e(f_2,g,r))\yv & (\text{by } (\ref{trivial_lower_bound})) \\
   & \geq \lambda_{\min}(C) - l_+ + \yv^{\top} (-g N_j) \yv \\
   & = \lambda_{\min}(C) - l_+ + 2 g (t-x_e) x_e \\
   & \geq \lambda_{\min}(C) - l_+ + 2 g \epsilon^2_0 t_0^2 & (\text{by } \ref{eq:case3_bounded_last}) \\
   & \geq 0 & (\text{by } \ref{rule:5}).
\end{align*}
This completes the proof that $U$ is copositive.

{Finally, let us verify the objective value of $U$. Since $\KK_i$ contributes to the objective value by $0$, each $G_j$ contributes to the objective value by $-r$, $H$ contributes to the objective value by $-(1+\sum_{i \in [m]} b_i^2)$ and $T$ contributes to the objective value by $-1$, the direct calculation leads to that the objective value of $U$ is $l- t_2 \epsilon_0 \cdot (1+\sum\limits_{i \in [m]} b_i^2) - r\cdot|\mathcal{B}|$.}
\end{proof}

We now consider the case where the feasible set is unbounded and $Q$ is PSD.
The construction is similar to the one in Theorem~\ref{thm_cop_structure}, but the argument is slightly more complex
since Lemmas~\ref{lemma_bounded_strict} and~\ref{lem_lower_bound_t_bounded_case} fail in the unbounded case.
Instead, we will make use of Theorem \ref{thm_final_robust} and Lemma \ref{lem_small_copostive block}, which hold when $Q$ is PSD.
We will also use that $Q$ can be written as $VV^{\top}$ for some $V$.

\begin{theorem} (Part~\ref{p2} of Theorem~\ref{thm:allstrong} - unbounded case)\label{thm_cop_structure_unbound}
Let $l$ be the lower bound of the optimal value that $l \leq z(\bv)$. Let $t_1,t_2$ be the same {as} in Theorem \ref{thm_final_robust} and let $u_j(\cdot)$ the same {as} in Remark \ref{rem_u_j}.
%When $\P$ is unbounded but 
When $Q$ is PSD, there exists some $\rho,\tau > 0$ that only depend on $A,\bv,\cv,Q,\mathcal{B}$, such that
for any $\epsilon_0 \in (0,t_1)$ and any $r > 0$,
\begin{align*}
   U := U(f_1, f_2, g, r, \alpha, l) = C + f_1 \left(\sum\limits_{i \in [m]} \KK_i\right) + \sum\limits_{j \in \mathcal{B}} G_j(f_2,g,r) + \tau H - l T 
\end{align*}
is copositive for $f_1, f_2, g, \tau$ satisfying the following rules:
\begin{enumerate}[label=(\text{rule}.\roman*),align=left]
    \item $f_2,g,r$ satisfies the condition in Lemma  \ref{lem_small_copostive block} so that $G_j(f_2,g,r)$ is copositive, \label{rule:unbounded_1}
    \item  $\tau = t_2 \epsilon_0$, \label{rule:unbounded_2}
    \item  $f_1 \geq \max \left\{  \rho + l, \frac{1}{2 \tau}, \frac{\rho+l}{\epsilon^2_0}, \max\limits_{j \in \mathcal{B}}\left\{ \frac{\rho+l}{u        _j(\epsilon_0)^2}  \right\}  \right\},$ \label{rule:unbounded_3}
    \item $g \geq \frac{\rho + l}{\epsilon^2_0}.$ \label{rule:unbounded_4}
\end{enumerate}
Moreover, {for the given $\epsilon_0,r,l$, the objective value of $U$ is $l- t_2 \epsilon_0 \cdot (1+\sum\limits_{i \in [m]} b_i^2) - r\cdot|\mathcal{B}|$.} This objective value can be arbitrarily {close} to $l$ as $r$ and $\epsilon_0$ goes to zero.
\end{theorem}

\begin{proof}
Let $\yv := [t ; \xv] \geq 0$,
since Lemma \ref{lem_small_copostive block} still holds by rule.i, we can still express $\yv^{\top} U \yv$ in the following way:
{
\begin{align*}
   \yv^{\top} U \yv & = \yv^{\top} C \yv + \sum_{i \in [m]} 
 \underbrace{\yv^{\top} f_1 (\KK_i) \yv}_{ \geq 0} + \sum_{j \in \mathcal{B}} \underbrace{\yv^{\top} ( G_j(f_2,g,r))\yv}_{\geq 0} + \underbrace{\yv^{\top}\tau H \yv}_{\geq 0} - \yv^{\top} l T \yv
\end{align*}
}
Again, the only term that is potentially negative is $\yv^{\top} (-lT) \yv + \yv^{\top} C \yv = -l t^2 + 2t \cv^{\top} \xv + \xv^{\top} Q  
 \xv$. Since $Q$ is PSD, this term is non-negative when $t = 0$. Therefore when $t = 0$, $\yv^{\top} U \yv \geq 0$. Therefore, we may assume that $t = 1$.
We will consider two cases:
\begin{enumerate}[label=(\text{case}.\arabic*),align=left]
    \item $[1 ; \xv] \in \P(\epsilon_0) \cap \I(\epsilon_0)$.
    \item $[1 ; \xv] \notin \P(\epsilon_0) \cap \I(\epsilon_0)$,
\end{enumerate}
where $\P(\epsilon_0)$ and $\I(\epsilon_0)$ are defined in (\ref{eq:idefn}) and (\ref{eq:pdefn}).

(case.1)
If $[1;\xv] \in \P(\epsilon_0) \cap \I(\epsilon_0)$, then $\xv$ is a feasible solution of MBQP($\epsilon_0$). Since $\epsilon_0 < t_1$ by our choice, Theorem \ref{thm_final_robust} implies that
\begin{align}
 2\cv^{\top} \xv + \xv^{\top} Q \xv \geq l - t_2 \epsilon_0 \label{eq:case_1_QPSD}.
\end{align}
Thus it follows that \begin{align*}
\yv^{\top} U \yv & \geq \yv^{\top} C \yv + \yv^{\top} (-lT) \yv +  \yv^{\top} \tau H \yv \\
& = 2\cv^{\top} \xv + \xv^{\top} Q \xv +  (-l) + \yv^{\top} \tau H \yv \\
& \geq  2\cv^{\top} \xv + \xv^{\top} Q \xv + (-l) + \yv^{\top} \tau T \yv  \\
& \geq  l - t_2 \epsilon_0  - l + \tau & (\text{by } (\ref{eq:case_1_QPSD})) \\
& = -t_2 \epsilon_0 + \tau \\
& = 0 & (\text{by } \ref{rule:unbounded_2}).
\end{align*}
(case.2) If $[1;\xv] \not\in \P(\epsilon_0) \cap \I(\epsilon_0)$, the difficulty here is that $\yv^{\top} C \yv = 2 \cv^{\top} \xv + \xv^{\top} Q \xv$ could potentially go to negative infinity. This is different from the proof of Theorem \ref{thm_cop_structure} because $\yv$ is normalized in a different way.

Since $Q$ is PSD, there exists some $V \in \R^{r \times n}$ such that $Q = V^{\top} V$ and $\xv^{\top} Q \xv = (V \xv)^{\top}(V \xv)$. 
Consider the following linear program:
\begin{align}
    \begin{array}{rl}\label{LP_large_obj_larg_violated}
    & \tau := \min\limits_{\xv,\varphi} \varphi \\
    & \text{s.t. } 2 \cv^{\top} \xv  = -1, \\
    & \ \ \ \ \ -\varphi \ev \leq V \xv \leq \varphi \ev, \\
    & \ \ \ \ \ - \varphi \leq \av_i^{\top} \xv \leq \varphi,\forall i \in [m], \\
    & \ \ \ \ \ \xv \geq 0, \varphi \geq 0
    \end{array}
\end{align}

First observe that $\tau > 0$. If $\tau = 0$, this implies that there exists some $\dv$ such that $2 \cv^{\top} \dv < 0$  and  $V \dv = 0$ and $\av_i^{\top} \dv = 0,\forall i \in [m]$ and $\dv \geq 0$,  implying that  $d_j = 0,\forall j \in \mathcal{B}$. Pick any feasible solution $\xv^*$ of (MBQP), one can verify that $\xv^*  + t \dv$ remains feasible for all $t \geq 0$ and its objective value goes to negative infinity as $t$ goes to infinity. This shows that 
the optimal value of the original (MBQP) is unbounded, which leads to a contradiction.
Since (\ref{LP_large_obj_larg_violated})is feasible and bounded from below, we have that $\tau > 0$ exists.

Select $\lambda_0$ such that
%quadratic $\tau^2 \lambda^2 - 2\lambda  - l$ such that 
for all $\lambda \geq \lambda_0$, $\tau^2 \lambda^2 - 2\lambda  - l \geq 0$ and then select
\begin{align}
\label{eq:choice_rho}
  \rho := \max\{\lambda_0,\max\limits_{i \in [m]}\{ \frac{|b_i| +1}{\tau} \}\}.
\end{align}
 % We will consider two subcases:
 % \begin{enumerate}[label=(\text{subcase}.\arabic*),align=left]
 %     \item $2 \cv^{\top} \xv \geq -\rho.$
 %     \item $2 \cv^{\top} \xv < -\rho$.
 % \end{enumerate}
If $2 \cv^{\top} \xv \geq -\rho$, since $\yv  \notin \P(\epsilon_0) \cap \I(\epsilon_0)$, applying the same argument in the proof of Theorem \ref{thm_cop_structure}, one of the following must hold:
\begin{align*}
    & \exists e_1 \in [m], |(\av_{e_1}^{\top} \xv)  - b_{e_1}| > \epsilon_0,  \\
    & \exists e_2 \in [m],j_2 \in \mathcal{B}, x_{j_2}  > 1 + \epsilon_0 \textup{implying that }  |(\av_{e_2}^{\top} \xv)  - b_{e_2}| > u_{j_2}(\epsilon_0), \\
    & \exists j_3 \in \mathcal{B}, x_{j_3} \in (\epsilon_0,1-\epsilon_0).
\end{align*}
This implies one of the following must hold:
\begin{align*}
    & \exists e_1 \in [m], \yv^{\top} (f_1 \KK_{e_1}) \yv \geq \epsilon_0^2 f_1,  \\
    & \exists e_2 \in [m],j_2 \in \mathcal{B}, \yv^{\top} (f_1 \KK_{e_2}) \yv \geq (u_{j_2}(\epsilon_0))^2 f_1, \\
    & \exists j_3 \in \mathcal{B}, \yv^{\top} (G_{j_3}(f_2,g,r)) \yv \geq \yv^{\top} N_{j_3} \yv \geq 2 \epsilon_0^2 g
\end{align*}
Since $\yv^{\top} (-lT) \yv + \yv^{\top} C \yv \geq -\rho - l$, with (rule.iii) and (rule.iv), one can assert that $\yv^{\top} U \yv \geq 0$.

It remains to consider the case when $2 \cv^{\top} \xv < -\rho$.
We may write $\xv = \lambda \xv_1$ where $2 \cv^{\top} \xv_1 = -1$ and $\lambda > \rho$. Since $\xv_1$ is a feasible solution in (\ref{LP_large_obj_larg_violated}) and by definition $\tau$, one of the following must hold:
 % \begin{enumerate}[label=(\text{subcase2}.\alph*),align=left]
 %     \item $\exists j \in [n], (V \xv_1)_{j} \geq \tau$ \label{subcase2a}.
 %     \item $\exists i \in [m], \av_i^{\top} \xv_1 \geq \tau$ \label{subcase2b}.
 % \end{enumerate}
\begin{align}
    & \exists j_4 \in [n], |(V \xv_1)_{j_4}| = \tau, \label{eq:suba} \\
    & \exists e_5 \in [m], |\av_{e_5}^{\top} \xv_1| = \tau \label{eq:subb}
\end{align}
If (\ref{eq:suba}) occurs, since $\lambda > \rho \geq \lambda_0$ by (\ref{eq:choice_rho}), it follows that
     \begin{align*}
         \yv^{\top} U \yv & \geq \yv^{\top} (-lT) \yv + \yv^{\top} C \yv  \\
          & = 2 \cv^{\top} \xv + \xv^{\top} Q \xv -l \\
          & = 2 \lambda \cv^{\top}  \xv_1 + \lambda^2 (\xv_1)^{\top} Q \xv_1 -l \\
          & \geq  -2\lambda + \lambda^2 \tau^2 -l & (\text{by (\ref{eq:suba})})\\
          & \geq 0 & (\text{by our choice of } \lambda_0 \text{ and } \lambda > \lambda_0)
     \end{align*}
If (\ref{eq:subb}) occurs, since $2 \cv^{\top} \xv < -\rho$, let $2 \cv^{\top} \xv = -\rho - \xi$. Then we may write $\xv = \xv_2 + \xi \xv_1$ where $\xv_1,\xv_2 \geq 0,\xi >0,\xv_2 = \rho \xv_1$. 
     Since $\rho \geq \max\limits_{i \in [m]}\{ \frac{|b_i| + 1}{\tau}\}$ by (\ref{eq:choice_rho}), then it follows that
     \begin{align}
        \label{unbounded_last}
       |\av^{\top}_{e_5} \xv - b_{e_5}| = |\av^{\top}_{e_5} \xv_2 + \xi \av^{\top}_{e_5} \xv_1 - b_{e_5}| \geq 1 + \tau \xi. 
     \end{align}
This implies 
     \begin{align*}
          \yv^{\top} U \yv & \geq \yv^{\top} (-lT) \yv + \yv^{\top} C \yv + \yv^{\top} f_1 \KK_{e_5} \yv  \\
          & \geq (-\rho -l - \xi) + f_1  (1+ \xi \tau)^2 & \text{(by (\ref{unbounded_last}))}\\
          & = f_1 \tau^2 \xi^2 + (2 f_1 \tau-1) \xi + f_1 -\rho -l \\
          & \geq 0 & \text{(by \ref{rule:unbounded_3})}
     \end{align*}
     The last inequality comes from the fact that this lower bound $f_1 \tau^2 \xi^2 + (2 f_1 \tau-1) \xi + f_1 -\rho -l$ is a quadratic function on $\xi$ and as long as $f_1 -\rho -l \geq 0, 2 f_1 \tau-1 \geq 0$, we can ensure $\yv^{\top} U \yv \geq 0$ for all $\xi \geq 0$ and this condition is implied by  \ref{rule:unbounded_3}.
     
{Finally, let us verify the objective value of $U$. Since $\KK_i$ contributes to the objective value by $0$, each $G_j$ contributes to the objective value by $-r$, $H$ contributes to the objective value by $-(1+\sum_{i \in [m]} b_i^2)$ and $T$ contributes to the objective value by $-1$, the direct calculation leads to that the objective value of $U$ is $l- t_2 \epsilon_0 \cdot (1+\sum\limits_{i \in [m]} b_i^2) - r\cdot|\mathcal{B}|$.}
\end{proof}

{
\begin{remark}
    Both Theorem \ref{thm_cop_structure} and Theorem \ref{thm_cop_structure_unbound} construct a feasible solution $U$ in (\ref{COP}) with objective value $l- t_2 \epsilon_0 \cdot (1+\sum\limits_{i \in [m]} b_i^2) - r\cdot|\mathcal{B}|$ for any given $\epsilon_0 \in (0,t_1)$ and $r > 0$. To construct a dual feasible solution promised by
    Theorem~\ref{thm:allstrong}\ref{p2} with objective value $l-\epsilon$ for some $\epsilon > 0$, one can use the construction of Theorem \ref{thm_cop_structure} and Theorem \ref{thm_cop_structure_unbound} by choosing $\epsilon_0 \in (0,t_1)$ and $r > 0$ such that
    \begin{align*}
        t_2 \epsilon_0 \cdot (1+\sum\limits_{i \in [m]} b_i^2) + r\cdot|\mathcal{B}| = \epsilon.
    \end{align*}
\end{remark}
}
\subsection{Theorem~\ref{thm:allstrong}\ref{p4}: (COP-dual) is not attainable in general.}~\label{sec:notatt}

Consider the maximum stable set problem for a graph $G = (V,E)$.
%Clearly, the optimal solution is just a single vertex.
This can be written as the following (MBQP):
\begin{align*}
    & \min -2\sum_{j \in V} x_i \\
    & \text{s.t } x_u + x_v + s_{e} = 1,\forall e := \{u,v\} \in E \\ 
    &  \ \ \ \ \ \xv \in \{0,1\}^{V}, \sv \geq 0
\end{align*}

By Theorem \ref{thm_cop_structure}, strong duality holds between its (CP-primal) and (COP-dual).
Its copositive dual (COP-dual) is 
{
\begin{align}
  \label{cop_dual_clique_standard}
  \begin{array}{rl}
   & \max - (\sum_{e \in E} 2 \alpha_e + \beta_e) - \theta \\
   & \text{s.t. } C + \sum_{e \in E} (\alpha_e A_e + \beta_e AA_e) + (\sum_{j \in V} \gamma_j N_j ) + \theta T = M \\
   & \ \ \ \ \ M \in \COP
   \end{array}
\end{align}
}
%Without losing generality, we may substitute $A_e$ with $\KK_e  :=  T  - A_e + AA_e$. In this case, we can write (COP-dual) as
{We can write the above (COP-dual) in the following equivalent form:}
{
\begin{align}
\label{cop_dual_clique}
    \begin{array}{rl}
 & \max - (\sum\limits_{e \in E}  \beta_e) - \theta \\
   & \text{s.t. } C + \sum\limits_{e \in E} (\mu_e \KK_e + \beta_e AA_e)  + (\sum\limits_{j \in V} \gamma_j N_j ) + \theta T = M \\
   & \ \ \ \ \ M \in \COP
    \end{array}
\end{align}
}

{
To see the equivalence between (\ref{cop_dual_clique_standard}) and (\ref{cop_dual_clique}), given any feasible solution $(\alpha^*,\beta^*,\gamma^*,\theta^*)$ in (\ref{cop_dual_clique_standard}), it straightforward to see that 
\begin{align*}
    & \mu_e = - \alpha^*_e,\forall e \in E, \\
    & \beta_e = \beta_e^* + \alpha^*_e, \forall e \in E, \\
    & \gamma_j = \gamma^*_j,\forall j \in V, \\
    & \theta = \theta^* + \sum_{e \in E} \alpha^*_e,
\end{align*}
is a feasible solution in (\ref{cop_dual_clique}) with the same objective. On the other hand, given any feasible solution $(\mu^*,\beta^*,\gamma^*,\theta^*)$ in (\ref{cop_dual_clique}), it straightforward to see that
\begin{align*}
    & \alpha_e = - \mu^*_e,\forall e \in E, \\
    & \beta_e = \beta_e^* + \mu_e^*, \forall e \in E, \\
    & \gamma_j = \gamma^*_j,\forall j \in V, \\
    & \theta = \theta^* + \sum_{e \in E} \mu^*_e,
\end{align*}
is a feasible solution in (\ref{cop_dual_clique_standard}) with the same objective.
}

\textbf{(Proof of Part~\ref{p4} of Theorem~\ref{thm:allstrong}):}
\begin{proof}
Consider the above COP problem for the special case where the graph $G$ is a clique of size six.
The stability number of a clique is one,
so the optimal value of \eqref{cop_dual_clique} is~$-2$.
Suppose the value of \eqref{cop_dual_clique} is attained.
Then there exists some $\mu_e^*,\beta^*_e,\theta^*,\gamma_j^*$ such that
\begin{align*}
    &  (\sum_{e \in E}  \beta_e^*) + \theta^* = 2 \\
    & M^* := M(\mu^*_e,\beta^*_e,\theta^*,\gamma_j^*) \in \COP
\end{align*}

For the sake of contradiction,
we will construct some non-negative vectors $\{\yv_i := [t_i;\xv_i;\sv_i] \geq 0\}$ such that $(\yv_i)^{\top} M^* \yv_i \geq 0$ {cannot} hold simultaneously.
We will construct those vectors sequentially. 
    We first show that {there exists some $k \in V$ such that $\gamma^*_k$ is significantly negative}. To see this,
    choose $\yv_1 := [1;\frac{1}{2}\ev;0 ]$,
    In this case, we have:
{
\begin{align*}
& \yv^{\top}_1 M^*\yv_1 = \yv^{\top}_1 \left(C + \sum\limits_{e \in E} \left(\mu_e^* \KK_e + \beta_e^* AA_e\right)  + \left(\sum\limits_{j \in V} \gamma_j^* N_j\right) + \theta^* T\right) \yv_1 \\
& = \underbrace{\yv^{\top}_1 C \yv_1}_{=-6} + \left(\sum\limits_{e \in E} \mu_e^* \underbrace{\yv^{\top}_1  \KK_e \yv_1}_{= 0} + \beta_e^* \underbrace{\yv^{\top}_1  AA_e \yv_1}_{=1}\right)  + \left(\sum\limits_{j \in V} \gamma_j^* \underbrace{\yv^{\top}_1N_j \yv_1}_{=-1/2}\right) + \theta^* \underbrace{\yv_1^{\top} T \yv_1}_{=1} \\
& = -6+\underbrace{\left(\sum_{e \in E}  \beta^*_e + \theta^*\right)}_{=2} -\frac{1}{2} \sum_{j \in V} \gamma^*_j \geq 0 \\
& \implies \sum_{j \in V} \gamma_j^* \leq -8 \\
& \implies \text{there exists some $\gamma^*_k \leq -\dfrac{4}{3}$} \text{ for some $k \in V$ since $|V| = 6$}.
\end{align*}
}
{The last inequality holds since if $\gamma^*_j > -\dfrac{4}{3}$ for all $j$, then we would have $\sum_{j \in V} \gamma_j^* > -8 $.}    
% \begin{align*}
% & \yv^{\top}_1 M^*\yv_1 = -6+\underbrace{(\sum_{e \in E}  \beta^*_e + \theta^*)}_{=2} -\frac{1}{2} \sum_{j \in V} \gamma^*_j \geq 0 \\
% & \implies \sum_{j \in V} \gamma_j^* \leq -8 \\
% & \implies \text{there exists some $\gamma^*_k \leq -\dfrac{4}{3}$} \text{ {for some $k \in V$}}.
% \end{align*}
% Without losing generality, we may assume that $\gamma^*_1$ is the smallest among $\gamma_j^*$ and $\gamma^*_1 \leq -\dfrac{4}{3}$.

Now for sufficiently small $\epsilon > 0$ and pick any arbitrary $v \in V$, we  
choose $\yv_{+}$ in the following way:
\begin{align*}
    t = 1+ \epsilon, x_j = \begin{cases}
        1 & \text{ if } j = v \\
        0 & \text{ otherwise}
    \end{cases}
    , s_e = \begin{cases}
        0 & \text{ if } e \in \delta(v) \\
        1 & \text{ otherwise},
    \end{cases}.
\end{align*}
Since $M^*$ is copositive, we have $(\yv_{+})^{\top} M^* \yv_{+} \geq 0$
\begin{align*}
   \implies  & -2(1+\epsilon) + ( (1+\epsilon)^2 \theta^* + \left(\sum_{e \in E}  \beta^*_e)\right) - 2\gamma^*_v \epsilon + \underbrace{O(\epsilon^2)}_{\text{introduced by } \yv_{+}^{\top} \KK_e \yv_{+}} \geq 0 \\
   \implies  & -2(1+\epsilon) + (2\epsilon \theta^* + \epsilon^2 \theta^*)  + \underbrace{( \theta^* + \left(\sum_{e \in E}  \beta^*_e\right) }_{=2} - 2\gamma^*_v \epsilon + O(\epsilon^2) \geq 0 \\
   \implies & -2 \epsilon + 2\epsilon \theta^* - 2\gamma^*_v \epsilon + {O(\epsilon^2)} \geq 0 \\
   \implies & -1 +  \theta^* - \gamma^*_v  + O(\epsilon) \geq 0 .
\end{align*}
Applying the same idea, we choose $\yv_{-}$ by {replacing} $\epsilon$ with $-\epsilon$
and then derive $$-1 +  \theta^* - \gamma^*_v  + O(\epsilon) \leq 0.$$ Combining with previous result, we get
\begin{align}
\label{eq:linear_not_att}
& -O(\epsilon) \leq \theta^* - \gamma^*_v - 1 \leq O(\epsilon) 
\implies \theta^* - \gamma^*_v - 1 = 0,\forall v \in V.
\end{align}
This means all $\gamma^*_v$ are the same. Since $\gamma_k^* \leq -\dfrac{4}{3}$, this implies that $\theta^* < 0$ and further implies that
\begin{align*}
    \sum_{e \in E} \beta_e^* = 2 - \theta^* > 0
\end{align*}
Given $i \in V$, and let {$\delta(i) := \{j : (i,j) \in E\}$} and 
define $G_i := \sum_{e \in \delta(i)} \beta_{{e}}^*$. Summing up all $G_i$ for $i \in V$, it yields that
\begin{align*}
    \sum_{i \in V} G_i = 2 \sum_{e \in E} \beta_e^* = 2 (2-\theta^*)
\end{align*}
Since $|V| = 6$ { and similar to the argument above}, this implies that there exists some $b \in V$ such that $G_b \leq \dfrac{1}{3} (2-\theta^*)$. Without losing generality, we may assume that $b = 1$.
For the last vector, for sufficiently small $\epsilon > 0$, we construct $\yv_2$ in the following way:
\begin{align*}
    t = 1, x_j = \begin{cases}
        1+\epsilon & \text{ if } j = 1 \\
        0 & \text{ otherwise}
    \end{cases}, s_e = \begin{cases}
        0 & \text{ if } e \in \delta(1) \\
        1 & \text{ otherwise }
    \end{cases}.
\end{align*}
Since $M^*$ is copositive, we have $\yv_2^{\top} M^* \yv_2 \geq 0, \\$
\begin{align*}
     \implies & -2(1+\epsilon) +  \theta^* + \left(\sum_{e \in \delta(1)}  (1+\epsilon)^2 \beta^*_e\right) + \left(\sum_{e \notin \delta(1)}  \beta^*_e\right)  + 2 \gamma^*_1(1+\epsilon) \epsilon + \underbrace{O(\epsilon^2)}_{\text{introduced by } (\yv_2)^{\top} \KK_e \yv_2} \geq 0 \\
\implies & -2(1+\epsilon) +  \underbrace{ \theta^* + \left(\sum_{e \in E}  \beta^*_e\right) }_{=2}  
 +  \underbrace{\left(\sum_{e \in \delta(1)} 2 \epsilon  \beta^*_e\right) }_{= \frac{2}{3} \epsilon (2-\theta^*) } + \underbrace{\left(\sum_{e \in \delta(1)}  \epsilon^2  \beta^*_e\right)   }_{= O(\epsilon^2)}   
  + 2 \gamma^*(1+\epsilon) \epsilon + O(\epsilon^2) \geq 0 \\
     \implies & -2 \epsilon + \frac{2}{3} \epsilon (2-\theta^*) + 2\gamma_1^* (1+\epsilon) \epsilon + O(\epsilon^2) \geq 0 \\
     \implies & -2  + \frac{2}{3}  (2-\theta^*) + 2\gamma_1^* (1+\epsilon)  + O(\epsilon) \geq 0 \\
      \implies & -2  + \frac{2}{3}  (1-\gamma_1^*) + 2\gamma_1^* (1+\epsilon)  + O(\epsilon) \geq 0 \ \ \text{(by (\ref{eq:linear_not_att}))} \\
      \implies & -2  + \frac{2}{3}  (1-\gamma_1^*) + 2\gamma_1^*   + O(\epsilon) \geq 0 \\
      \implies & -\frac{4}{3} + \frac{4}{3} \gamma_1^* + O(\epsilon) \geq 0 \ \text{ (contradiction since $\gamma_1^* < 0$)}.
     % \implies & -\theta + 2 \gamma_i(1+\epsilon) + O(\epsilon^2) \geq 0 \\
     % \implies & -1 - \gamma_i + 2 \gamma_i(1+\epsilon) + O(\epsilon^2) \gtrapprox 0 \\
     % \implies & -1 + \gamma_i(1+\epsilon) + \gamma_i \epsilon + O(\epsilon^2) \gtrapprox 0 \\
     % \implies & \text{ contradiction }
\end{align*}

\end{proof}

The above result implies that when $\P$ is unbounded and $Q$ is PSD then (\ref{CP}) may not have a Slater point in the general case. 
This can also be shown through a simpler construction:
let $Q = ([n-1;-\ev]) ([n-1;-\ev])^{\top} \in \SS^{n}$, consider the following instance of (MBQP):
\begin{align*}
\min \{\xv^{\top} Q \xv | \xv \geq 0\}. \end{align*}
By (Part~\ref{p2} of Theorem~\ref{thm:allstrong}), strong duality holds. Its (COP-dual) takes the form of
\begin{align}
    \label{dual_no_slater}
    & \max - t \notag \\
    & \text{s.t. } \begin{bmatrix}
        t & 0 \\
        0 & Q
    \end{bmatrix} = M, M \in \COP.
\end{align}
One can see that there is no Slater point in (\ref{dual_no_slater}) since $[0;1;\ev]^{\top} M [0;1;\ev] = 0$ no matter what $t$ is.

\section{Conclusion and future directions.}\label{sec:conclude}
We proved sufficient conditions for strong duality to hold between Burer's reformulation of MBQPs and its dual. 
{By replacing the copositive constraint by $\SS_+ + \SS_P$, we have proposed an SDP-based algorithm to conduct sensitivity analysis of general (MBQP), which provides much better bounds than existing methods.} This algorithm is motivated by the structure of $\epsilon$-optimal solution of the COP dual. However, the sizes of instances we can currently perform sensitivity analysis are limited by the SDP solver.  

One direction of research is to extend such strong duality results for reformulations of more general QCQPs~\cite{burer2012representing}. 
{Since we consider replacing the copositive constraint  by $\SS_+ + \SS_P$ an interesting open question is to provide some theoretical guarantee on such modification~\cite{klep2024random}.}
%, we cannot perform sensitivity analysis on practice instances and have to illustrate our idea with small random instances. 
Another possible future direction is to develop a more scalable solver for the SDP in (\ref{eq:actual1}),
for instance, using the techniques from \cite{majumdar2020recent}. 
{Finally, recent works consider solving (\ref{COP}) using cutting-plane techniques~\cite{ANSTREICHER2021218,badenbroek2022analytic,guocopositive,LinderothRaghunathan}. 
It may be possible that the structure information provided by part~\ref{p2} of Theorem~\ref{thm:allstrong} could offer opportunities to design more efficient cutting-plane techniques.}

\section*{Acknowledgements}
We would like to thank the support by Office of Naval Research grant N000142212632 and the National Science Foundation (NSF) grant 2112533 (AI4OPT).
\bibliography{ref}

\appendix

\section{Dual bounds are hard to find}
\label{sec:dual_bounds_hard_to_find}
{In Table~\ref{tab:review}, we provide numerical evidence that “modern solvers are able to find a good solution of (MPQP) in a reasonable time while finding it more challenging to provide matching dual bounds.” We use the instance from \cite{cifuentes2024lagrangian} and we generate ten random instances. For time $T \in \{60s,1200s\}$, we report the ratio of the best primal value found by Gurobi in $T$s over the best known dual bound (first column) and the ratio of the best primal value found by Gurobi in $T$s over the best dual bound found by Gurobi in $T$s (second column). In all ten instances, Gurobi finds a nearly optimal solution in the 60s, and there is no change in the best primal solution afterward. On the other hand, it is very hard for Gurobi to provide matching dual bounds, and the dual gap is more than 10\% in the 1200s.}

{\footnotesize
\begin{table}[h!]\centering
\caption{avg performance of PATH-STAB}
\label{tab:review}
% \ra{0.5}
\setlength{\tabcolsep}{3pt}
\begin{tabular}{@{}c|cccc@{}}\toprule
time   & $\frac{\text{best primal value}}{\text{best known dual bound}}$ &  $\frac{\text{best primal value}}{\text{best dual bound of Gurobi}}$ \\
\midrule
60s & 1.6\%    &  17.2\%  \\
1200s & 1.6\%    &  12.1\%  \\
\bottomrule
\end{tabular}
\end{table}
}

\section{Details of computational results}\label{Appdix:comp}

This appendix provides more refined information regarding the tables from section~\ref{sec:comp_body}.
Specifically,
Tables~\ref{tab:comb1}--\ref{tab:comb3} expand on 
Table~\ref{tab:comb} by providing results for each density level.
Similarly,
Tables~\ref{tab:sslpd03}--\ref{tab:sslpd07} expand on 
Table~\ref{tab:lpind}, Tables~\ref{tab:qpind1}--\ref{tab:qpind3} expand on 
Table~\ref{tab:qpind}, and Tables~\ref{tab:PACK1}--\ref{tab:PACK3} expand on 
Table~\ref{tab:PACK}.

{\footnotesize
\begin{table}[!htbp]
\centering
\ra{0.5}
\setlength{\tabcolsep}{4pt}
\caption{Average relative gap (COMB) -- density=0.3}
\label{tab:comb1}
\begin{tabular}{@{}cccccccccccccc@{}}\toprule 
$\Delta k$ & 1 & 2 & 3 & 4 & 5 & 6 & 7 & 8 & 9 & 10 & \phantom{a} & avg time(s)\\
\midrule
  Shor1 & 1  & 1
  & 1
  & 1
  & 1
  & 1
  & 1
  & 1
  & 1
  & 1
  & &2.80 \\
 Shor2   & 2.13	& 4.97&	6.64	&9.30	&12.24&	8.26	&6.81&	6.25	&5.87&	5.75  && 3.88 \\
 Our method   &  1.35	& 0.14& 0.04	& 0.00	& 0.07	&0.26	&0.35	&0.45	&0.50	& 0.54 && 4.03\\
  Our method ($\sigma =2$)  & 0.99&	0.00&		0.43&		1.10&		2.06	&	1.87	&	1.84&		1.91&		1.95&		2.02 && 3.63\\
  
    Our method ($\sigma =4$)  & 5.26 &	1.22	& 0.20 &	0.00	& 0.19 & 	0.48 &	0.63 &	0.79 &	0.89 &	0.98 && 3.89\\
        Our method ($\sigma =6$)  & 19.97 &	7.47	&2.75 &	1.05	 & 0.18 &	0.00 &	0.05 &	0.17	 &0.27&	0.36 && 4.01\\
         Our method ($\sigma =8$)  & 35.30&15.82&	6.46	&3.28&	1.45	&0.29&	0.03&	0.00	&0.03	&0.06 && 3.48\\
          Our method ($\sigma =10$)  & 45.71	&21.99	&10.18	&6.01	&3.44	&0.95	&0.29	&0.09	&0.02 & 0.00 && 3.91\\
 Cont   &  0.80 &
 0.93 &
 1.03 &
 1.02 &
 1.05 &
 1.04 &
 1.03 &
 1.03 &
 1.03 &
 1.02 && 0.00 \\
\bottomrule
\end{tabular}
\end{table}
}

{\footnotesize
\begin{table}[!htbp]
\centering
\ra{0.5}
\setlength{\tabcolsep}{4pt}
\caption{Average relative gap (COMB) -- density=0.5}
\label{tab:comb2}
\begin{tabular}{@{}cccccccccccccc@{}}\toprule 
$\Delta k$ & 1 & 2 & 3 & 4 & 5 & 6 & 7 & 8 & 9 & 10 & \phantom{a} & avg time(s)\\
\midrule
  Shor1 & 1  & 1
  & 1
  & 1
  & 1
  & 1
  & 1
  & 1
  & 1
  & 1
  & &2.98 \\
 Shor2   & 1.50 &	3.14 &	5.02	&5.56&	6.08	&6.26&	6.19	&6.21&	6.28	&6.22  && 4.45 \\
 our method   &  1.13	&0.38&	0.02&	0.00	&0.03&	0.15&	0.23&	0.30&	0.36&	0.41 && 4.77\\
  Our method ($\sigma =2$)  & 0.54	& 0.023	& 0.25	& 0.75 & 	1.20	& 1.56 & 	1.80	& 2.00 & 	2.18	& 2.29 && 4.17\\
  
    Our method ($\sigma =4$)  & 4.40 	&1.13&	0.16&	0.00	&0.09&	0.31	&0.52&	0.69&	0.85&	0.98 && 5.93\\
        Our method ($\sigma =6$)  & 11.44 &	3.73&1.55	&0.45	&0.06	&0.00&	0.053	&0.12	&0.21&	0.30 && 3.88\\
         Our method ($\sigma =8$)  & 17.58&6.72	&3.40	&1.32&	0.47	&0.13&	0.02&	0.00&	0.02&	0.07 && 4.83\\
          Our method ($\sigma= 10$)  & 28.91&11.74	&6.78	&3.14	&1.54	&0.71	&0.31	&0.11	&0.02&	0.00 && 4.77\\
 Cont   & 0.91 &
 1.05 &
 1.07 &
 1.06 &
 1.05 &
 1.04 &
 1.02 &
 1.03 &
 1.03 &
 1.02 && 0.00 \\
\bottomrule
\end{tabular}
\end{table}
}

{\footnotesize
\begin{table}[!htbp]
\centering
\ra{0.5}
\setlength{\tabcolsep}{4pt}
\caption{Average relative gap (COMB) -- density=0.7}
\label{tab:comb3}
\begin{tabular}{@{}cccccccccccccc@{}}\toprule 
$\Delta k$ & 1 & 2 & 3 & 4 & 5 & 6 & 7 & 8 & 9 & 10 & \phantom{a} & avg time(s)\\
\midrule
  Shor1 & 1  & 1
  & 1
  & 1
  & 1
  & 1
  & 1
  & 1
  & 1
  & 1
  & &3.63 \\
 Shor2   & 1.27 &	3.02&	4.95&	5.86&	5.95&	5.75&	5.75&	5.78&	5.67&	5.56  && 4.21 \\
 Our method   &  1.36 & 	0.53&	0.12&	0.00&	0.11&	0.23&	0.30&	0.37&	0.44&	0.50 && 5.01\\
 Our method ($\sigma =2$)  & 0.59	&0.03&0.30&0.85&1.37&1.67&	1.93&	2.15&2.28&	2.39 && 6.28\\
  
    Our method ($\sigma =4$)  & 4.64 &	1.34 &	0.27	&0.00&	0.18&	0.41&	0.61&	0.78&	0.94	&1.00 && 5.16\\
        Our method ($\sigma =6$)  & 11.56&	4.62&	2.03	&0.49&	0.06&	0.00	&0.04&	0.12&	0.23&	0.34 && 4.82\\
         Our method ($\sigma =8$)  & 18.91	& 7.88	& 4.18	& 1.68	& 0.61& 	0.19& 	0.04& 	0.00& 	0.03& 	0.10 && 6.12\\
          Our method ($\sigma =10$)  & 30.84& 	13.89	&8.25&	3.91	&1.76&	0.75&	0.30&	0.08&	0.01&	0.00 && 4.54\\
 Cont   &  1.00 &
 1.07 &
 1.09 &
 1.02 &
 1.02 &
 1.06 &
 1.05 &
 1.05 &
 1.04 &
 1.04 && 0.00 \\
\bottomrule
\end{tabular}
\end{table}
}

{\footnotesize
\begin{table}[!htbp]
\setlength{\tabcolsep}{6pt}
\centering
\caption{Average relative gap (SSLP) -- density=0.3}
\label{tab:sslpd03}
 \begin{tabular}{ccccc} 
 \toprule
$\norm{\Delta \bv}_{\infty}$ &$\leq 1$ & $\leq 2$ & $\leq 3$ & avg time(s)  \\  
  \midrule
 Shor1 & 1 & 1 & 1 & 3.68 \\
 Shor2 &  1.25   & 1.63 &  1.83 & 7.14 \\
 our method  & 0.54  & 0.50  &  0.56 &  5.64 \\
 Cont & 1.00 & 1.00 & 1.00 & 0.00 \\ [1ex] 
 \hline
 \end{tabular}
\end{table}
}

{\footnotesize
\begin{table}[!htbp]
\setlength{\tabcolsep}{6pt}
\centering
\caption{Average relative gap (SSLP) -- density=0.5}
\label{tab:sslpd05}
 \begin{tabular}{ccccc} 
 \toprule
$\norm{\Delta \bv}_{\infty}$ &$\leq 1$ & $\leq 2$ & $\leq 3$ & avg time(s)  \\  
  \midrule
 Shor1 & 1 & 1 & 1 & 3.53 \\
 Shor2 &1.09 &  1.18 &  1.26 & 7.11 \\
 Our method   & 0.68  & 0.64  & 0.68 &  5.68 \\
 Cont & 1.00 & 1.00 & 1.00 & 0.00 \\ [1ex] 
 \hline
 \end{tabular}
\end{table}
}

{\footnotesize
\begin{table}[!htbp]
\setlength{\tabcolsep}{6pt}
\centering
\caption{Average relative gap (SSLP) -- density=0.7}
\label{tab:sslpd07}
 \begin{tabular}{ccccc} 
 \toprule
$\norm{\Delta \bv}_{\infty}$ &$\leq 1$ & $\leq 2$ & $\leq 3$ & avg time(s)  \\  
  \midrule
 Shor1 & 1 & 1 & 1 & 3.69 \\
     Shor2 &  1.03   & 1.08 &  1.11 & 7.38 \\
 Our method    & 0.74  & 0.70 & 0.73 &   6.15\\
 Cont & 1.00 & 1.00 & 1.00 & 0.00 \\ [1ex] 
 \hline
 \end{tabular}
\end{table}
}

{\footnotesize
\begin{table}[!htbp]
\setlength{\tabcolsep}{6pt}
\centering
\caption{Average relative gap (SSQP) -- density=0.3}
\label{tab:qpind1}
 \begin{tabular}{ccccc} 
 \toprule
$\norm{\Delta \bv}_{\infty}$ &$\leq 1$ & $\leq 2$ & $\leq 3$ & avg time(s)  \\  
  \midrule
 Shor1 & 1 & 1 & 1 & 3.53\\
 Shor2 &  1.35  &  1.56 &   1.76 &  7.11 \\
 Our method  & 0.52   & 0.48  &  0.53 &  6.12 \\
 Cont & 1.00 & 1.00 & 1.00 & 0.00 \\ [1ex] 
 \hline
 \end{tabular}
\end{table}
}

{\footnotesize
\begin{table}[!htbp]
\setlength{\tabcolsep}{6pt}
\centering
\caption{Average relative gap (SSQP) -- density=0.5}
% \label{tab:qpind2}
 \begin{tabular}{ccccc} 
 \toprule
$\norm{\Delta \bv}_{\infty}$ &$\leq 1$ & $\leq 2$ & $\leq 3$ & avg time(s)  \\  
  \midrule
 Shor1 & 1 & 1 & 1 &  3.59\\
 Shor2 &  1.02 &  1.06 &  1.1 &  7.10 \\
 Our method  &  0.51  & 0.48 &  0.52 &  5.56\\
 Cont & 1.00 & 1.00 & 1.00 & 0.00 \\ [1ex] 
 \hline
 \end{tabular}
\end{table}
}

{\footnotesize
\begin{table}[!htbp]
\setlength{\tabcolsep}{6pt}
\centering
\caption{Average relative gap (SSQP) -- density=0.7}
\label{tab:qpind3}
 \begin{tabular}{ccccc} 
 \toprule
$\norm{\Delta \bv}_{\infty}$ &$\leq 1$ & $\leq 2$ & $\leq 3$ & avg time(s)  \\  
  \midrule
 Shor1 & 1 & 1 & 1 &  3.62\\
 Shor2 &  0.97  & 0.98 &  1.0 &  7.00 \\
 Our method  & 0.63  & 0.62 &  0.67 &   6.02 \\
 Cont & 1.00 & 1.00 & 1.00 & 0.00 \\ [1ex] 
 \hline
 \end{tabular}
\end{table}
}

\begin{table}[!htbp]
    {\footnotesize
\caption{Average relative gap of 2000 samples for (PACK)-- density=0.3}
\label{tab:PACK1}
% \ra{0.5}
\setlength{\tabcolsep}{20pt}
\centering
 \begin{tabular}{ccc} 
 \toprule
$\norm{\Delta \bv}_{\infty}$ & $\leq 3$ & time(s)  \\  
  \midrule
 Shor1 & 1  &908 \\
 Shor2 &  1.30   &1199 \\
 Our method  0.61 & 0.60     &	 1007 \\
 Cont & 1.00 &  0 \\ [1ex] 
 \bottomrule
 \end{tabular}

}
\end{table}

\begin{table}[!htbp]
    {\footnotesize
\caption{Average relative gap of 2000 samples for (PACK)-- density=0.5}
\label{tab:PACK2}
% \ra{0.5}
\setlength{\tabcolsep}{20pt}
\centering
 \begin{tabular}{ccc} 
 \toprule
$\norm{\Delta \bv}_{\infty}$ & $\leq 3$ & time(s)  \\  
  \midrule
 Shor1 & 1  &892 \\
 Shor2 &  1.53   &1061 \\
 Our method  0.62 & 0.60     &	 1028 \\
 Cont & 1.00 &  0 \\ [1ex] 
 \bottomrule
 \end{tabular}

}
\end{table}

\begin{table}[!htbp]
    {\footnotesize
\caption{Average relative gap of 2000 samples for (PACK)-- density=0.7}
\label{tab:PACK3}
% \ra{0.5}
\setlength{\tabcolsep}{20pt}
\centering
 \begin{tabular}{ccc} 
 \toprule
$\norm{\Delta \bv}_{\infty}$ & $\leq 3$ & time(s)  \\  
  \midrule
 Shor1 & 1 &876 \\
 Shor2 &  7.78  &1172 \\
 Our method  0.58 & 0.60     &	 970 \\
 Cont & 0.96 &  0 \\ [1ex] 
 \bottomrule
 \end{tabular}

}
\end{table}

\end{document}